\documentclass[aap,MSNbibl,dvips]{arximspdf}
\usepackage{graphics}
\doi{10.1214/09-AAP652}
\volume{20}
\issue{3}
\pubyear{2010}
\firstpage{1068}
\lastpage{1097}

\makeatletter
\newproclaim{definition}{Definition}[section]
\newtheorem{theorem}[definition]{Theorem}
\newtheorem{proposition}[definition]{Proposition}
\newtheorem{corollary}[definition]{Corollary}
\newtheorem{lemma}[definition]{Lemma}
\newtheorem{remark}[definition]{Remark}
\makeatother

\begin{document}
\begin{frontmatter}

\title{Topology-guided sampling of nonhomogeneous random processes}
\runtitle{Topology-guided random field sampling}

\begin{aug}
\author[A]{\fnms{Konstantin} \snm{Mischaikow}\thanksref{a1}\ead[label=e1]{mischaik@math.rutgers.edu}} \and
\author[B]{\fnms{Thomas} \snm{Wanner}\thanksref{a2}\ead[label=e2]{wanner@math.gmu.edu}\corref{}}
\runauthor{K. Mischaikow and T. Wanner}
\affiliation{Rutgers University and George Mason University}
\address[A]{Department of Mathematics \\
Rutgers University \\
Piscataway, New Jersey 08854\\ USA \\
\printead{e1}} 
\address[B]{Department of Mathematical Sciences \\
George Mason University \\
Fairfax, Virginia 22030\\ USA \\
\printead{e2}}
\thankstext{a1}{Supported in part by NSF Grants
DMS-05-11115 and DMS-01-07396, by DARPA and by the U.S.~Department
of Energy.}
\thankstext{a2}{Supported in part by NSF Grant DMS-04-06231
and the U.S.~Department of Energy under
Contract DE-FG02-05ER25712.}
\end{aug}

\received{\smonth{4} \syear{2008}}
\revised{\smonth{10} \syear{2009}}

%
\begin{abstract}
Topological measurements are increasingly being accepted as an important
tool for quantifying complex structures. In many applications, these
structures can be expressed as nodal domains of real-valued functions
and are obtained only through experimental observation or numerical simulations.
In both cases, the data on which the topological measurements are based are
derived via some form of finite sampling or discretization. In this
paper, we
present a probabilistic approach to quantifying the number of
components of
generalized nodal domains of nonhomogeneous random processes on the real
line via finite discretizations, that is, we consider excursion sets of a
random process relative to a nonconstant deterministic threshold function.
Our results furnish explicit probabilistic a priori bounds for the suitability
of certain discretization sizes and also provide information for the
choice of
location of the sampling points in order to minimize the error
probability. We illustrate our results for a variety of random processes,
demonstrate how they can be used to sample the classical nodal domains
of deterministic functions perturbed by additive noise and discuss
their relation to the density of zeros.
\end{abstract}

%

\end{frontmatter}
%

\section{Introduction}
\label{secintro}
The~motivation for this work comes from our attempts to create novel metrics
for quantifying, comparing and cataloging large sets of complicated varying
geometric patterns. Random fields (for a general background,
see~\cite{adler81a,bharuchareids86a,farahmand98a,kahane85a,marcuspisier81a}, as well as the references therein)
provide a framework in which to approach these problems and
have, over the last few decades, emerged as an important
tool for studying spatial phenomena which involve an element of
randomness~\cite{adler81a,adlertaylor07a,rozanov98a,torquato02b,vanmarcke83a}. For the types of applications, we
have in mind~\cite{gameiroetal04a,gameiroetal05a,krishanetal07a},
we are often satisfied with a topological classification of sub- or
super-level sets
of a scalar function. Algebraic topology, and in particular homology,
can be used
in a computationally efficient manner~\cite{kaczynskietal04a} to
coarsely quantify
these geometric properties. In past work~\cite{dayetal07a,mischaikowwanner07a},
we developed a probabilistic framework for assessing the correctness of homology
computations for random fields via uniform discretizations. The~approach considers
the homology of nodal domains of random fields which are given by classical
Fourier series in one and two space dimensions, and it provides
explicit and
sharp error bounds as a function of the discretization size and
averaged Sobolev
norms of the random field. While we do not claim it is trivial---there are
complicated combinatorial questions that need to be resolved---we believe
that it is possible to extend the methods and hence the results
of~\cite{mischaikowwanner07a} to higher-dimensional domains.

\begin{figure}

\includegraphics{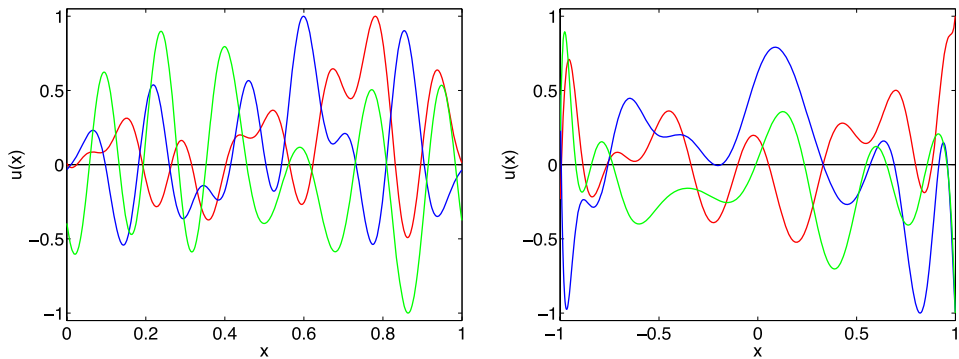}

\caption{Sample functions from a random sum of the form
$u(x,\omega) = \sum_{k=0}^N g_k(\omega) \varphi_k(x)$
where~$g_1, \ldots, g_N$ are independent standard Gaussian
random variables. In the left diagram, we consider random periodic
functions, that is, the basis functions~$\varphi_k$ are given by
$\varphi_{2k}(x) =\break  \cos(2\pi k x)$ and $\varphi_{2k-1}(x) = \sin
(2\pi k x)$,
in the right diagram they are the Chebyshev polynomials
$\varphi_k(x) = \cos(k \arccos x)$. In each case, we choose $N = 16$.}
\label{figsample}
\end{figure}

The~more serious restriction in~\cite{mischaikowwanner07a}
is the use of periodic random fields, which due to the fact that the associated
spatial correlation function is homogeneous, simplifies many of the estimates.
In general, however, one expects to encounter nonhomogeneous random fields.
In such cases, it seems unreasonable to expect that uniform sampling
provides the optimal choice. For example, in Figure~\ref{figsample},
three sample
functions each are shown for a random sum involving periodic basis
functions and
Chebyshev polynomials. As one would expect, the zeros of the random Chebyshev
sum are more closely spaced at the boundary, and therefore small uniform
discretization are most likely not optimal for determining the topology
of the
nodal domains.

With this as motivation, we allow for a more general sampling
technique. We
remark that because of the subtlety of some of the necessary estimates we
restrict our attention in this paper to one-dimensional domains.

\begin{definition}[(Nonuniform approximation of generalized nodal domains)]
\label{defcubappr}
Consider a compact interval $[a,b]\subset\mathbb{R}$, a
threshold function
$\mu\dvtx  [a,b] \to\mathbb{R}$, and a function $u \dvtx  [a,b] \to\mathbb
{R}$. Then we
define the \textit{generalized nodal domains of~$u$} by
%
\begin{equation} \label{defnodal}
N_\mu^\pm= \bigl\{ x \in[a,b]   \dvtx    \pm\bigl( u(x) - \mu(x)
\bigr)
\ge0 \bigr\},
\end{equation}
which for the case of $\mu(x) \equiv0$ reduces to the classical definition
of a nodal domain in~\cite{couranthilbert53a}.
An \textit{$M$-discretization of~$[a,b]$} is a collection of~$M + 1$
grid points
\[
a = x_0 < x_1 < \cdots< x_M = b   ,
\]
and we define~$x_{M+1} = x_M = b$ in the following.
The~\textit{cubical approximations~$Q_M^\pm$} of the \textit{generalized
nodal domains~$N_\mu^\pm$ of~$u$} are defined as the sets
\[
Q_{\mu,M}^\pm:= \bigcup\bigl\{ [x_k,x_{k+1}]   \dvtx
\pm\bigl( (u - \mu)(x_{k}) \bigr) \ge0   ,
k = 0,\ldots,M \bigr\} .
\]
\end{definition}

Given a subset $X\subset[a,b]$, let~$\beta_0(X)$ denote the number of
components
of~$X$. Consider a random field $u \dvtx  [a,b] \times\Omega\to\mathbb
{R}$ over the
probability space $(\Omega, \mathcal{F}, \mathbb{P})$. We are
interested in optimally
characterizing the topology, that is, determining the number of
components, of
the nodal domains $N_\mu^\pm$ in terms of the cubical\vspace*{-1.5pt}
approximations~$Q_{\mu,M}^\pm$. In other words, our goal is to choose the
$M$-discretization of~$[a,b]$ in such a way as to optimize
\[
\mathbb{P}\{ \beta_0(N_\mu^\pm)=\beta_0(Q_{\mu,M}^\pm
)\}
  .
\]
We provide two results addressing this question. The~first
characterizes the choice of the sampling points $a = x_0 < x_1 < \cdots
< x_M = b$
under reasonably general abstract conditions. More precisely, consider the
following assumptions:
\begin{enumerate}[(A3)]
\item[(A1)] For every $x \in[a,b]$, we have $\mathbb{P}\{ u(x) = \mu
(x) \}
= 0$.
\item[(A2)] The~random field is such that $\mathbb{P}\{ u - \mu
\mbox{ has
a double zero in }  [a,b] \} = 0$.
\item[(A3)] For $\sigma\in\{ \pm1 \}$, $x \in[a,b]$ and $\delta>
0$ with
$x + \delta\in[a,b]$ define
\begin{eqnarray*}
p_\sigma(x,\delta) &=& \mathbb{P}\biggl\{
\sigma u(x) \ge\sigma\mu(x) ,
\sigma u\biggl(x + \frac{\delta}{2}\biggr) \le\sigma
\mu\biggl(x + \frac{\delta}{2}\biggr) ,\\
&&\hspace*{107pt}\sigma u(x + \delta) \ge\sigma\mu(x + \delta) \biggr\} .
\end{eqnarray*}
Then there exists a continuously differentiable function
$\mathcal{C}_0 \dvtx  [a,b] \to\mathbb{R}^+$ as well as a constant
$\mathcal{C}_1 > 0$ such
that for all $x \in[a,b]$ with $x + \delta\in[a,b]$ we have
\[
p_{+1} (x,\delta) + p_{-1} (x,\delta)   \le
\mathcal{C}_0(x) \cdot\delta^3 + \mathcal{C}_1 \cdot\delta^4
  .
\]
\end{enumerate}

In Section~\ref{seclocal}, we prove the following result.

\begin{theorem}[(Sampling based on local probabilities)]
\label{thmTopSamp}
Consider a probability space $(\Omega, \mathcal{F}, \mathbb{P})$,
a continuous threshold function $\mu\dvtx  [a,b] \to\mathbb{R}$,
and a random field
$u \dvtx  [a,b] \times\Omega\to\mathbb{R}$ over~$(\Omega, \mathcal{F},
\mathbb{P})$ such that
for $\mathbb{P}$-almost all $\omega\in\Omega$ the function
$u(\cdot,\omega) \dvtx  [a,b] \to\mathbb{R}$ is continuous.
Choose the sampling points~$a = x_0 < \cdots< x_M = b$
such that
\[
\int_{x_{k-1}}^{x_k} \sqrt[3]{\mathcal{C}_0(x)}  \, dx   =
\frac{1}{M} \cdot\int_a^b \sqrt[3]{\mathcal{C}_0(x)} \,  dx
\qquad  \mbox{for all }
k = 1,\ldots,M   ,
\]
and consider the generalized nodal domains~$N_\mu^\pm(\omega)$ and their
approximations~$Q_{\mu,M}^\pm(\omega)$ as in Definition~\textup{\ref{defcubappr}}.
If assumptions~\textup{(A1)}, \textup{(A2)} and~\textup{(A3)} hold, then
%
\begin{equation} \label{thmTopSamp1}
\mathbb{P}\{ \beta_0(N_\mu^\pm) = \beta_0(Q_{\mu,M}^\pm)
\} \geq
1 - \frac{4}{3 M^2} \cdot\biggl( \int_a^b \sqrt[3]{\mathcal
{C}_0(x)} \,  dx
\biggr)^3
+ O\biggl( \frac{1}{M^3} \biggr)   .
\end{equation}
\end{theorem}

This theorem is a direct generalization of the corresponding result
in~(\cite{mischaikowwanner07a}, Theorem~1.3). Numerical computations presented
in Section~\ref{secappl} suggest that for certain nonhomogeneous random
fields this estimate is sharp---and in fact an enormous improvement over
the homogeneous result where~$\mathcal{C}_0(x)$ is replaced by~$\max
_{x \in G} \mathcal{C}_0(x)$.

Of course in practice one is interested in applying Theorem~\ref{thmTopSamp}
to specific random fields. This requires the verification of assumptions~(A1),
(A2) and~(A3), preferably in terms of central random field characteristics.

\begin{definition}
For a random field~$u \dvtx  [a,b] \times\Omega\to\mathbb{R}$ over
a probability
space~$(\Omega, \mathcal{F}, \mathbb{P})$, we define its \textit{spatial correlation
function}~$R \dvtx  [a,b]^2 \to\mathbb{R}$ as
\[
R(x,y) = \mathbb{E}\bigl( \bigl( u(x) - \mathbb{E}u(x) \bigr)
\bigl( u(y) - \mathbb{E}u(y) \bigr) \bigr)\qquad
 \mbox{for all }
x,y \in[a,b] ,
\]
where~$\mathbb{E}$ denotes the expected value of a random variable
over~$(\Omega, \mathcal{F}, \mathbb{P})$.
\end{definition}

If the random field is sufficiently smooth, then the derivatives of the
spatial correlation function,
%
\begin{equation} \label{assRFcorr1}
R_{k,\ell}(x) =
\frac{\partial^{k + \ell} R}{\partial x^k \,\partial y^\ell}(x,x),
\end{equation}
have a natural interpretation in terms of spatial derivatives of the
random field~$u$. Since
\[
R_{k,\ell}(x)    =
{\mathbb{E}} \bigl( \bigl( u^{(k)}(x) - \mathbb{E}u^{(k)}(x) \bigr)
\bigl( u^{(\ell)}(y) - \mathbb{E}u^{(\ell)}(y) \bigr) \bigr) ,
\]
the function~$R_{k,k}$ contains averaged information on the square of the
$k$th derivative of the random function~$u$, more precisely, its variance.

To relate the spatial correlation function to the function~$\mathcal
{C}_0$ in
Theorem~\ref{thmTopSamp}, we specialize to Gaussian random fields. To be
more precise, we make the following assumptions.
\begin{enumerate}[(G2)]
\item[(G1)] Consider a Gaussian random field $u \dvtx  [a,b] \times\Omega
\to\mathbb{R}$
over a probability space~$(\Omega, \mathcal{F}, \mathbb{P})$ such
that $u(\cdot
,\omega) \dvtx
[a,b] \to\mathbb{R}$ is twice continuously differentiable for
$\mathbb{P}$-almost all
$\omega\in\Omega$. Furthermore, assume that for every $x \in[a,b]$ the
expected value of~$u(x)$ satisfies
\[
\mathbb{E}u(x) = 0 .
\]
\item[(G2)] The~spatial correlation function~$R$ is three times continuously
differentiable in a neighborhood of the diagonal~$x = y$ and the matrix
%
\begin{equation} \label{assRFcorr2}
\mathcal{R}(x) =
\pmatrix{
R_{0,0}(x) & R_{1,0}(x) & R_{2,0}(x)\vspace*{2pt} \cr
R_{1,0}(x) & R_{1,1}(x) & R_{2,1}(x)\vspace*{2pt} \cr
R_{2,0}(x) & R_{2,1}(x) & R_{2,2}(x)
}
\end{equation}
is positive definite for all $x \in[a,b]$.
\end{enumerate}

We make considerable use of~$\mathcal{R}$, and thus introduce the following
notation
%
\begin{eqnarray}
\mathcal{R}_{33}^m & := & R_{0,0} R_{1,1} - R_{1,0}^2   , \nonumber\\
\mathcal{R}_{32}^m & := & R_{0,0} R_{2,1} - R_{1,0} R_{2,0}   ,
\label
{defRFcorr}\\
\mathcal{R}_{31}^m & := & R_{1,0} R_{2,1} - R_{1,1} R_{2,0}   .
\nonumber
\end{eqnarray}
These expressions are just the determinants of minors of~$\mathcal{R}$.
This allows us to state the following theorem.
\begin{theorem}[(Sampling based on spatial correlation)]
\label{thmCorrSamp}
Consider a Gaussian random field~$u \dvtx  [a,b] \times\Omega\to\mathbb{R}$
satisfying~\textup{(G1)} and~\textup{(G2)}, and a threshold function $\mu\dvtx  [a,b] \to
\mathbb{R}
$ of
class~$C^3$. Choose the sampling points~$a = x_0 < \cdots< x_M = b$
in such a way that
\[
\int_{x_{k-1}}^{x_k} \sqrt[3]{\mathcal{C}(x)}\,   dx   =
\frac{1}{M} \cdot\int_a^b \sqrt[3]{\mathcal{C}(x)} \,  dx
\qquad  \mbox{for all }
k = 1,\ldots,M   ,
\]
where
%
\begin{equation} \label{thmCorrSamp1}
\mathcal{C}(x)   =   \frac{\det\mathcal{R}(x)}{48 \pi
\mathcal{R}_{3,3}^m(x)^{3/2}} \cdot
\bigl( 1 + \mathcal{A}(x) \bigr) \cdot e^{-\mathcal{B}(x)}   ,
\end{equation}
given
\begin{eqnarray*}
\mathcal{A}(x) & = & \frac{( \mathcal{R}_{3,1}^m(x) \mu(x) -
\mathcal{R}_{3,2}^m(x)
\mu^\prime(x)
+ \mathcal{R}_{3,3}^m(x) \mu^{\prime\prime}(x) )^2}{\mathcal
{R}_{3,3}^m(x)
\det\mathcal{R}(x)}   \ge  0   , \\
\mathcal{B}(x) & = & \frac{( R_{1,0}(x) \mu(x) - R_{0,0}(x)
\mu
^\prime(x) )^2 +
\mathcal{R}_{3,3}^m(x) \mu(x)^2}{2   R_{0,0}(x)   \mathcal
{R}_{3,3}^m(x)}   \ge
  0   .
\end{eqnarray*}
Let~$Q_{\mu,M}^\pm(\omega)$ denote the cubical approximations of the random
generalized nodal domains~$N_\mu^\pm(\omega)$ of~$u(\cdot,\omega
)$. Then
%
\begin{equation} \label{thmCorrSamp2}
\mathbb{P}\{ \beta_0(N_\mu^\pm)= \beta_0(Q_{\mu,M}^\pm)
\}
 \ge
1 - \frac{1}{M^2} \cdot\biggl(
\int_a^b \sqrt[3]{\mathcal{C}(x)}  \, dx \biggr)^3 +
O\biggl( \frac{1}{M^{3}} \biggr)   .
\end{equation}
\end{theorem}

The~proof of Theorem~\ref{thmCorrSamp} is presented in Section~\ref
{secspatial}.
However, it depends on nontrivial results concerning the asymptotic
behavior of
sign-distribution probabilities of parameter-dependent Gaussian random
variables.
These results are developed in Section~\ref{secPTool}.

The~number of nodal domains $\beta_0(N_\mu^\pm)$ is clearly
dependent upon
the zeros of~$u-\mu$. Thus, it is reasonable to expect that there is some
relationship between the function~$\mathcal{C}$ derived in
Theorem~\ref{thmCorrSamp}
and the density of the zeros of the random field~$u$. The~first step is to
obtain a density function. For this, a weaker form of~(G2) is sufficient.
\begin{enumerate}[(G3)]
\item[(G3)] Assume that the spatial correlation function~$R$ is two times
continuously differentiable in a neighborhood of the diagonal~$x = y$ and
that $R(x,x) > 0$ for all $x \in[a,b]$.
\end{enumerate}
Finding the density of the zeros of random fields has been studied in a variety
of settings, see, for example,~\cite{adlertaylor07a,bharuchareids86a,cramerleadbetter04a,edelmankostlan95a,farahmand98a}, as well as the
references therein. The~following theorem can be found
in~\cite{cramerleadbetter04a}, (13.2.1), page~285.

\begin{theorem}[(Density of zeros of a random field)]
\label{thmDensZero}
Consider a Gaussian random field~$u \dvtx  [a,b] \times\Omega\to\mathbb{R}$
satisfying~\textup{(G1)} and~\textup{(G3)}. Then the density function for the number of
zeros of~$u$ is given by
%
\begin{equation} \label{thmDensZero1}
\mathcal{D}(x) = \frac{\mathcal{R}_{3,3}^m(x)^{1/2}}{\pi\cdot
R_{0,0}(x)} .
\end{equation}
In other words, for every interval~$I \subset[a,b]$ the expected
number of zeros of~$u$
in~$I$ is given by $\int_I \mathcal{D}(x)  \, dx$.
\end{theorem}

While Theorem~\ref{thmDensZero} has been known for quite some time, its
implications are surprising. As is demonstrated through examples in
Section~\ref{secappl} there is no simple discernible relationship
between the function $\mathcal{C}^{1/3}$ of Theorem~\ref{thmCorrSamp}
and the
density function $\mathcal{D}$.

As is made clear at the beginning of this Introduction, our motivation
is to
develop optimal sampling methods for the analysis of complicated time-dependent
patterns. Thus, before turning to the proofs of the above-mentioned
results, we
begin, in Section~\ref{secappl}, with demonstrations of possible
applications and
implications of Theorem~\ref{thmCorrSamp}. In particular, we consider several
random generalized Fourier series $u \dvtx  [a,b] \times\Omega\to\mathbb{R}$
defined by
%
\begin{equation} \label{defRFourSer}
u(x,\omega) = \sum_{k=0}^\infty g_k(\omega) \cdot
\varphi_k(x)   ,
\end{equation}
where $\varphi_k \dvtx  [a,b] \to\mathbb{R}$, $k \in\mathbb{N}_0$,
denotes a family
of smooth functions and we assume that the Gaussian random
variables~$g_k \dvtx  \Omega\to\mathbb{R}$, $k \in\mathbb{N}_0$, are
defined over a
common probability space~$(\Omega, \mathcal{F}, \mathbb{P})$ with mean~$0$.

We conclude the paper with a general discussion of future work concerning
natural generalizations to higher dimensions.

\section{Sampling of specific random sums} \label{secappl}
To demonstrate the applicability and implications of
Theorem~\ref{thmCorrSamp}, we consider in this section several random
generalized Fourier series $u \dvtx  [a,b] \times\Omega\to\mathbb{R}$ of the
form in~(\ref{defRFourSer}). As mentioned before, the functions
$\varphi_k \dvtx  [a,b] \to\mathbb{R}$, $k \in\mathbb{N}_0$, denote a
family of smooth
functions and we assume that the random variables~$g_k \dvtx  \Omega\to
\mathbb{R}$,
$k \in\mathbb{N}_0$, are Gaussian with vanishing mean, and defined
over a common
probability space~$(\Omega, \mathcal{F}, \mathbb{P})$. We would like
to point out that
these random variables do not need to be independent, and we define
\[
\alpha_{k,m} = \mathbb{E}(g_k g_m)
\qquad  \mbox{for all }
k,m \in\mathbb{N}_0   .
\]
Then one can easily show that
\[
R_{k,\ell}(x)    =
{\mathbb{E}} \bigl( u^{(k)}(x) u^{(\ell)}(x) \bigr)    =
\sum_{i,j=0}^{\infty} \alpha_{i,j} \varphi_i^{(k)}(x) \varphi
_j^{(\ell)}(x)
  .
\]
If in addition the random variables~$g_k$ are pairwise independent,
then we have
\[
R_{k,\ell}(x)    =
\sum_{j=0}^{\infty} \alpha_{j,j} \varphi_j^{(k)}(x) \varphi
_j^{(\ell)}(x)
  ,
\]
where $\alpha_{j,j} \ge0$ for all $j \in\mathbb{N}_0$. One can show
that this
diagonalization can always be achieved for Gaussian random fields, provided
the basis functions~$\varphi_k$ are chosen appropriately. For more
details, we
refer the reader to~\cite{adlertaylor07a}, Theorems~3.1.1 and~3.1.2, Lemma~3.1.4.

Within the above framework of random generalized Fourier series, we
specifically consider several classes:
\begin{itemize}
\item\textit{Random Chebyshev polynomials} $u \dvtx  [-1,1] \times\Omega
\to\mathbb{R}$
of the form
%
\begin{equation} \label{exche}
u(x,\omega) = \sum_{k=0}^N g_k(\omega) \cdot
\cos( k \arccos x )
\qquad  \mbox{with }
\mathbb{E}( g_k g_\ell) = \delta_{k,\ell}   .
\end{equation}
\item\textit{Random cosine series} $u\dvtx  [0,1] \times\Omega\to\mathbb
{R}$ of the
form
%
\begin{equation} \label{excos}
u(x,\omega) = \sum_{k=0}^N g_k(\omega) \cdot
\cos( k \pi x )
 \qquad \mbox{with }
\mathbb{E}( g_k g_\ell) = \delta_{k,\ell}   .
\end{equation}
\item\textit{Random $L$-periodic functions} $u \dvtx  \mathbb{R}\times\Omega
\to
\mathbb{R}$
of the form
%
\begin{eqnarray} \label{exper}
u(x,\omega) &=& \sum_{k=0}^\infty a_k \cdot\biggl(
g_{2k}(\omega) \cdot\cos\frac{2\pi k x}{L} +
g_{2k-1}(\omega) \cdot\sin\frac{2\pi k x}{L} \biggr)\nonumber\\[-8pt]\\[-8pt]
\eqntext{\mbox{with }
\mathbb{E}( g_k g_\ell) = \delta_{k,\ell},}
\end{eqnarray}
with real constants~$a_k$.
\item\textit{Random polynomials} $u \dvtx  [-3,3] \times\Omega\to\mathbb{R}$
\textit{with Gaussian coefficients of binomial variance} of the form
%
\begin{equation} \label{expol}
u(x,\omega) = \sum_{k=0}^N g_k(\omega) \cdot x^k
\qquad  \mbox{with }
\mathbb{E}( g_k g_\ell) = \delta_{k,\ell} \cdot\pmatrix{ N
\cr k}
  ,
\end{equation}
\item\textit{Random polynomials} $u \dvtx  [-3,3] \times\Omega\to\mathbb{R}$
\textit{with Gaussian coefficients of unit variance} of the form
%
\begin{equation} \label{exalg}
u(x,\omega) = \sum_{k=0}^N g_k(\omega) \cdot x^k
\qquad  \mbox{with }
\mathbb{E}( g_k g_\ell) = \delta_{k,\ell}   .
\end{equation}
\end{itemize}
As is indicated in Section~\ref{secintro}, we assume that all the random
coefficients are centered Gaussian random variables over a common probability
space~$(\Omega,\mathcal{F},\mathbb{P})$.

\subsection{The~case of vanishing threshold function} \label{subsecmuzero}

We begin our applications by thresholding sample random sums at their expected
value, that is, we use the threshold function $\mu\equiv0$. In this particular
case, the function~$\mathcal{C}(x)$ defined by~(\ref{thmCorrSamp1}) in
Theorem~\ref{thmCorrSamp} simplifies to
%
\begin{equation} \label{Czero}
\mathcal{C}(x)   =   \frac{\det\mathcal{R}(x)}{48 \pi
\mathcal{R}_{3,3}^m(x)^{3/2}} ,
\end{equation}
since both~$\mathcal{A}(x)$ and~$\mathcal{B}(x)$ vanish.

\begin{figure}

\includegraphics{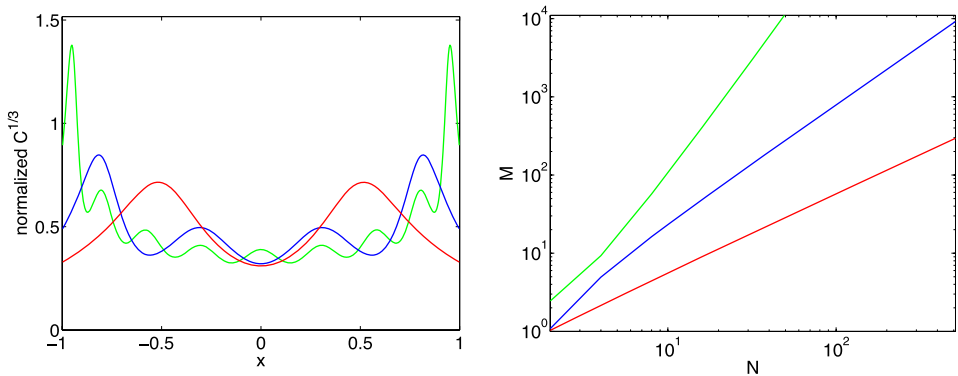}

\caption{Topology-guided sampling of random Chebyshev
polynomials~(\protect\ref{exche}).
The~left diagram shows the functions~$\mathcal{C}^{1/3}$ for $N = 3,
5, 10$
(red, blue and green, respectively---increasing values of~$N$
increase the number of extrema); for comparison reasons, each
curve has been scaled in such a way that the area under the graph
is one. The~right diagram shows the expected number of zeros of the
random Chebyshev polynomials as a function of~$N$ (bottom red curve),
the value of~$M$ for which Theorem~\protect\ref{thmCorrSamp} gives a correctness
probability of~$95\%$ (middle blue curve), and the value of~$M$ for
which~\protect\cite{mischaikowwanner07a} gives a correctness
probability of~$95\%$ (top green curve) with $\mathcal{C}_0 = \max
\mathcal{C}_0(x)$.}
\label{figcheb}
\end{figure}

For the case of random Chebyshev polynomials~(\ref{exche}), the left diagram
in Figure~\ref{figcheb} shows three normalized sample functions
\[
\frac{\mathcal{C}^{1/3}(x)}{\int_{-1}^1 \mathcal{C}^{1/3}(x) \, dx}
\]
for $N = 3,5,10$. The~right diagram shows the expected number of zeros
of the
random Chebyshev polynomials as a function of~$N$ (red curve), which grows
proportional to~$N$. Thus, in order to sample the random field sufficiently
fine, we expect to use significantly more than~$O(N)$ discretization points.
The~blue curve in the right diagram of Figure~\ref{figcheb} shows the
values of~$M$ for which the bound in~(\ref{thmCorrSamp2}) of
Theorem~\ref{thmCorrSamp} implies a correctness probability of~$95\%$,
and a
least squares fit of this curve furnishes $M \sim N^{3/2}$. For comparison,
the green curve in the same diagram shows the values of~$M$ for which
the bound
in our previous result~(\cite{mischaikowwanner07a}, Theorem~1.4)
implies a
correctness probability of~$95\%$, provided we apply this theorem
with~$\mathcal{C}_0$
given as the $\max_{x \in[-1,1]} \mathcal{C}_0(x)$. Notice that in
this case we
have~$M \sim N^3$. In other words, only the topology-guided sampling result
of the current paper yields a reasonable growth for the number of sampling
points. In fact, based on our results for periodic random fields
in~\cite{mischaikowwanner07a} and the numerical simulations
in~\cite{dayetal07a}, we expect that~$M \sim N^{3/2}$ is the optimal
discretization size.

\begin{figure}[b]

\includegraphics{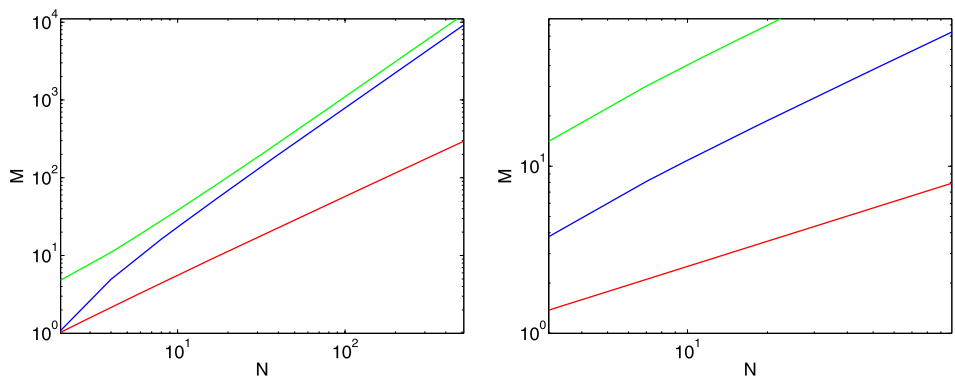}

\caption{Topology-guided sampling of random trigonometric
polynomials~(\protect\ref{excos}) satisfying Neumann boundary
conditions (left diagram) and random algebraic
polynomials~(\protect\ref{expol}) with binomial variances (right diagram).
The~curves show the expected numbers of zeros (bottom red curve),
the discretization size required by Theorem~\protect\ref{thmCorrSamp} to
achieve~$95\%$ correctness (middle blue curve), and the discretization
size required by~\protect\cite{mischaikowwanner07a} for a correctness
probability of~$95\%$ (top green curve), with $\mathcal{C}_0 = \max
\mathcal{C}_0(x)$.}
\label{figasymp}
\end{figure}

For the case of random cosine sums~(\ref{excos}), that is, random trigonometric
sums satisfying homogeneous Neumann boundary conditions, the analogue
of the
right diagram in Figure~\ref{figcheb} is depicted in the left diagram of
Figure~\ref{figasymp}. Notice that for the random cosine sums
the expected number of zeros is proportional to~$N$, and the required
number of sampling points has to be proportional to~$N^{3/2}$ for both
Theorem~\ref{thmCorrSamp} and~\cite{mischaikowwanner07a}, Theorem~1.4.
In other words, in this situation the gains from topology-guided
sampling are no longer as large as in the context of Chebyshev polynomials.
Also in this case, the curves for~$M$ are obtained in such a way that the
right-hand side in~(\ref{thmCorrSamp2}) or the corresponding bound
in~\cite{mischaikowwanner07a} equals~$95\%$

Similar behavior can be seen in the case of random polynomials~(\ref{expol})
with Gaussian coefficients of binomial variance; see the right diagram of
Figure~\ref{figasymp}. For the random algebraic polynomials~(\ref{expol}),
one can show that the expected number of zeros is proportional to~$N^{1/2}$,
and the required number of sampling points implied by~(\ref{thmCorrSamp2})
or~\cite{mischaikowwanner07a} has to be proportional to~$N^{3/4}$
for both
results. In fact, the function~$\mathcal{C}$ can be computed
explicitly in
this case.
Due to~(\ref{expol}), the spatial correlation function~$R$ is given by
\[
R(x,y)   =
\sum_{k=0}^N {N \choose k} x^k y^k   =
( 1 + xy )^N   ,
\]
which after some elementary computations furnishes
%
\begin{equation} \label{expol1}
\mathcal{C}(x)   =
\frac{N^{1/2} (N-1)}{24 \pi( 1 + x^2 )^3}   .
\end{equation}
As for the case of random polynomials with Gaussian coefficients of unit
variance, a classical result due to Kac~\cite{kac43a,kac49a}
implies that the expected number
of zeros is proportional to~$\log N$. In this case, Theorem~\ref{thmCorrSamp}
implies that the required number of sampling points has to be proportional
to~$( \log N )^{3/2}$.

\subsection{The~case of constant threshold function} \label{subsecmuconst}

We now turn our attention to a constant threshold function $\mu(x) =
\tau$,
for some real number~$\tau$. In this case, the function~$\mathcal
{C}(x)$ in
Theorem~\ref{thmCorrSamp} simplifies to
%
\begin{equation} \label{Cconst}
\mathcal{C}(x)   =   \frac{\det\mathcal{R}(x)}{48 \pi
\mathcal{R}_{3,3}^m(x)^{3/2}} \cdot\mathcal{S}(x)   ,
\end{equation}
where
%
\begin{equation} \label{Sconst}
\mathcal{S}(x)   =   \biggl( 1 + \frac{\mathcal{R}_{3,1}^m(x)^2
\tau^2}{\mathcal{R}
_{3,3}^m(x)
\det\mathcal{R}(x)} \biggr) \cdot
\exp\biggl( -\frac{R_{1,0}(x)^2 + \mathcal{R}_{3,3}^m(x)}{2
R_{0,0}(x)
\mathcal{R}_{3,3}^m(x)} \cdot\tau^2 \biggr)   .
\end{equation}

\begin{figure}

\includegraphics{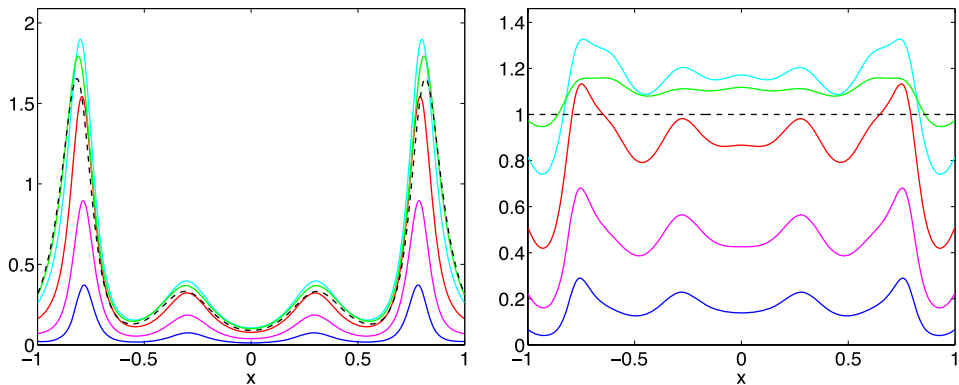}
\vspace*{-1mm}
\caption{Effect of varying the threshold~$\tau$ on the
function~$\mathcal{C}(x)$
in~(\protect\ref{Cconst}) for random Chebyshev polynomials~(\protect\ref{exche})
with~$N = 5$. The~left diagram shows the function~$\mathcal{C}(x)$ for
$\tau= 0, 1, 2, 3, 4, 5$ (black, green, cyan, red, magenta, blue),
the right diagram shows only the function~$\mathcal{S}(x)$ defined
in~(\protect\ref{Sconst}).}\vspace*{2mm}
\label{figchetau}
\end{figure}

\begin{figure}

\includegraphics{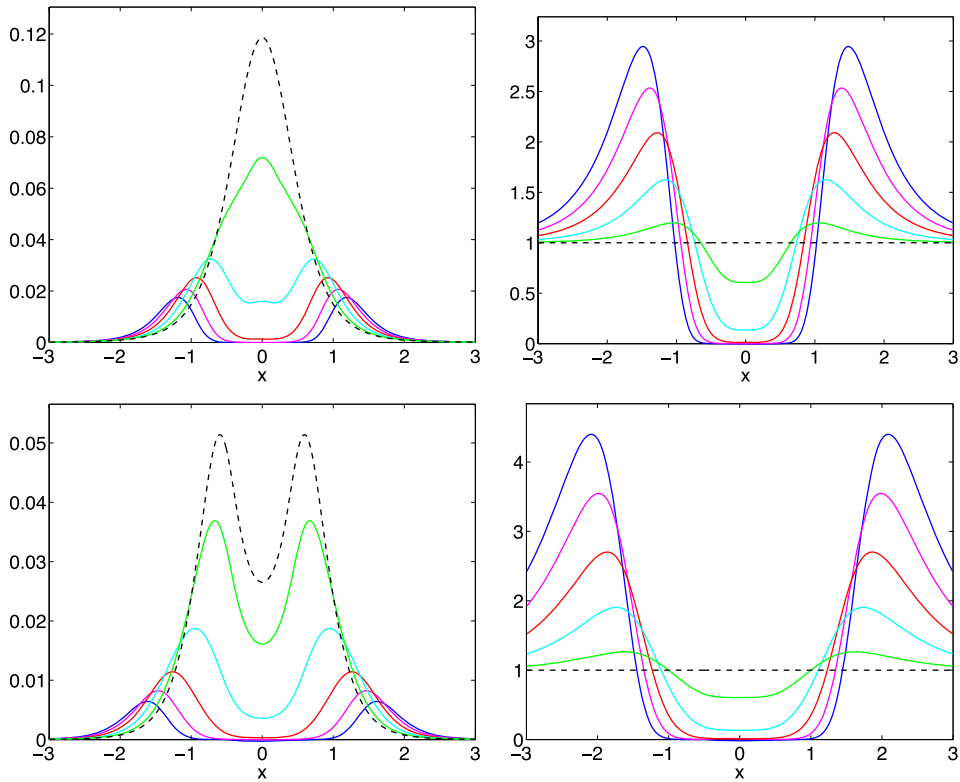}
\vspace*{-1mm}
\caption{Effect of varying the threshold~$\tau$ on the
function~$\mathcal{C}(x)$
in~(\protect\ref{Cconst}) for random algebraic polynomials~(\protect\ref{expol})
(top row) and~(\protect\ref{exalg}) (bottom row) with~$N = 5$. In each row,
the left diagram shows~$\mathcal{C}(x)$ for $\tau= 0, 1, 2, 3, 4, 5$ (black,
green, cyan, red, magenta, blue), and the right diagram shows only
the function~$\mathcal{S}(x)$ defined in~(\protect\ref{Sconst}).}
\label{figtaupolalg}
\end{figure}

For large values of~$|\tau|$, the scaling function~$\mathcal{S}(x)$
will be
close to zero,
and it therefore effectively decreases the probability for mistakes
in the homology computation. In fact, it decreases exponentially fast
with respect to~$|\tau|$. However, as is shown in Figure~\ref{figchetau}
for the random Chebyshev polynomials~(\ref{exche}), for values
of~$\tau$
close to zero, there can be regions in which the probability for mistakes
actually increases. This behavior is even more pronounced in the case of
random algebraic polynomials~(\ref{expol}) and~(\ref{exalg}), which is
shown in Figure~\ref{figtaupolalg}.

\subsection{The~case of varying threshold function} \label{subsecmuvary}

We now consider the case of a general
threshold function under the following assumptions.
Suppose a deterministic function~$\mu(x)$ is perturbed by a centered
Gaussian random field~$u(x,\omega)$, and that we are interested in
determining the classical nodal domains of the sum
\[
v(x,\omega) = \mu(x) + u(x,\omega)   .
\]
Sampling~$v(x,\omega)$ at the threshold zero is obviously equivalent
to sampling $u(x,\omega)$ at the threshold~$-\mu(x)$. Thus, we can use
Theorem~\ref{thmCorrSamp} to find the optimal location of the sampling
points using the function~$\mathcal{C}(x)$ defined in~(\ref{thmCorrSamp1}).

In order to demonstrate the effects of the varying threshold
function~$-\mu(x)$ more clearly, we now assume that the perturbing
random field~$u$ is homogeneous, that is, we assume that~$u$ is a random
$L$-periodic function of the form~(\ref{exper}). Furthermore, we assume
that the real scaling factors~$a_k$ in~(\ref{exper}) satisfy
\[
\sum_{k=0}^\infty k^6 a_k^2 < \infty  ,
\]
and that at least two of the~$a_k$ do not vanish. It was shown
in~\cite{mischaikowwanner07a} that in this case the spatial correlation
function~$R$ is given by
\[
R(x,y)  =
\mathbb{E}u(x) u(y)  =
\sum_{k=0}^\infty a_k^2 \cdot\cos\frac{2\pi k (x-y)}{L}   .
\]
From this, one can readily see that the matrix function~$\mathcal{R}(x)$
defined in~(\ref{assRFcorr2}) is constant and given by
\[
\mathcal{R}(x) =
\pmatrix{
A_0 & 0 & -\dfrac{4 \pi^2 A_1}{L^2} \vspace*{2pt}\cr
0 & \dfrac{4 \pi^2 A_1}{L^2} & 0 \vspace*{2pt}\cr
-\dfrac{4 \pi^2 A_1}{L^2} & 0 & \dfrac{16 \pi^4 A_2}{L^4}
},
\]
where
\[
A_\ell   =
\sum_{k=0}^{\infty} k^{2\ell} a_k^2   .
\]
Thus, the function~$\mathcal{C}(x)$ in~(\ref{thmCorrSamp1}) is now
given as
%
\begin{equation} \label{Cvaryper}
\mathcal{C}(x)   =   \frac{\pi^2}{6 L^3} \cdot
\frac{A_0 A_2 - A_1^2}{A_0^{3/2} A_1^{1/2}} \cdot\mathcal{S}(x)   ,
\end{equation}
where
\begin{eqnarray*}
\mathcal{S}(x)   &=&   \biggl( 1 + \frac{( A_1 \mu(x) + A_0
\mu^{\prime\prime}(x) \cdot({L^2}/{(4\pi^2)}) )^2}
{A_0 ( A_0 A_2 - A_1^2 )} \biggr) \\
&&{}\times
\exp\biggl( -\frac{A_1 \mu(x)^2 + A_0 \mu^\prime(x)^2 \cdot
({L^2}/{(4\pi^2)})}{2 A_0 A_1} \biggr)   .
\end{eqnarray*}
Notice that the exponential factor is bounded above
by~$\exp(-\mu(x)^2 / (2A_0))$, that is, large function values of~$\mu(x)$
lead to small failure probabilities.

\begin{figure}

\includegraphics{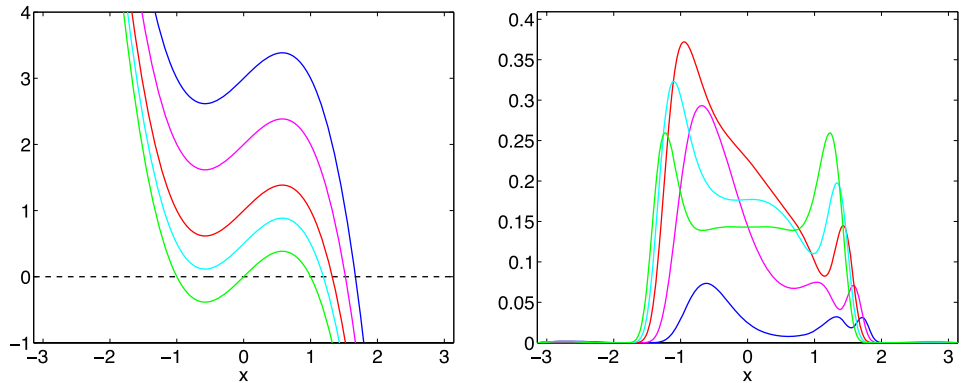}

\caption{Sampling of deterministic functions~$\mu(x)$ perturbed by
homogeneous random noise. The~left image shows the functions
$\mu(x) = x - x^3 + \tau$ for $\tau= 0, 0.5, 1, 2, 3$
(green, cyan, red, magenta, blue), the right images shows the
corresponding functions~$\mathcal{C}(x)$ defined in~(\protect\ref{Cvaryper}).}
\label{figpervary}
\end{figure}

We close this subsection by visualizing the function~$\mathcal{C}(x)$ defined
in~(\ref{Cvaryper}) for the deterministic function~$\mu(x) = x - x^3
+ \tau$
and $\tau$-values between~$0$ and~$3$. The~specific functions~$\mu
(x)$ are
shown in the left image of Figure~\ref{figpervary}. In the right image,
the corresponding functions~$\mathcal{C}(x)$ are shown, where~$u$ is
defined as
in~(\ref{exper}) with $a_k = 0$ for $k = 0$ and $k > N$, as well as
$a_k = N^{-1/2}$ for $k = 1,\ldots,N$. This implies that the variance
of~$u(x)$ equals~$1$. In Figure~\ref{figpervary}, we use $N = 5$.

\subsection{Comparison with density-guided sampling} \label{subsecdensity}

\begin{figure}

\includegraphics{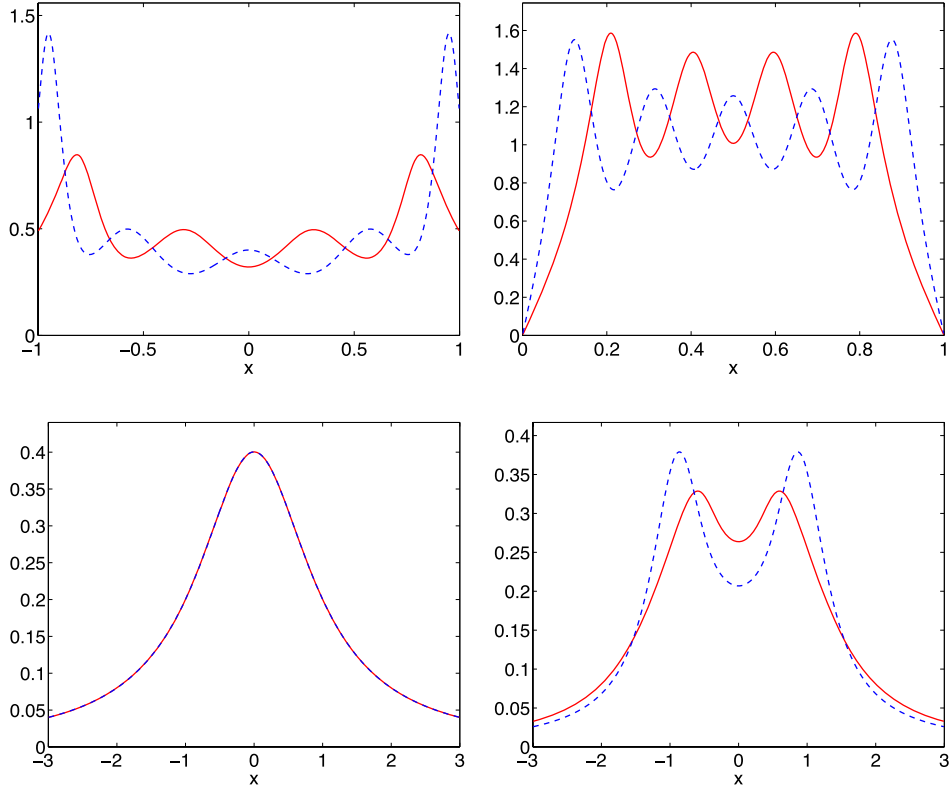}

\caption{A comparison of the function~$\mathcal{C}^{1/3}$ and the density
function~$\mathcal{D}$ for random Chebyshev polynomials~(\protect\ref{exche})
(top left diagram), random trigonometric polynomials~(\protect\ref{excos})
(top right), random algebraic polynomials~(\protect\ref{expol})
(bottom left), and random algebraic polynomials~(\protect\ref{exalg})
(bottom right). In all cases, the areas under the graphs have
been normalized to one, and we chose $N = 5$.}
\label{figdenshom}
\end{figure}

In order to illustrate the differences between the density of
zeros~$\mathcal{D}$
derived in Theorem~\ref{thmDensZero} and the function~$\mathcal
{C}^{1/3}$ from
Theorem~\ref{thmCorrSamp}, we return to our examples from the last section.
For each of these examples, Figure~\ref{figdenshom} depicts both
\[
\frac{\mathcal{C}^{1/3}(x)}{\int\mathcal{C}^{1/3}(x) \, dx}
\quad  \mbox{and}\quad
\frac{\mathcal{D}(x)}{\int\mathcal{D}(x) \, dx}
\]
for the case $N = 5$. It is evident
from these graphs that in most cases, the homology-based sampling
density is
different from the actual density of zeros. In fact, in many cases it behaves
anticyclic to~$\mathcal{D}$ in the sense that the local extrema
of~$\mathcal{C}^{1/3}$
alternate with the local extrema of~$\mathcal{D}$.

There are, however, exceptions, as the case of the random algebraic
polynomial~(\ref{expol}) demonstrates. In this case, it follows from
Theorem~\ref{thmDensZero} that
\[
\mathcal{D}(x)   =
\frac{N^{1/2}}{\pi( 1 + x^2 )}   ,
\]
and together with~(\ref{expol1}) this shows that the
normalized $\mathcal{C}^{1/3}$- and $\mathcal{D}$-functions coincide.

\section{Sampling based on local probabilities} \label{seclocal}
The~goal of this section is the proof of Theorem~\ref{thmTopSamp},
which is a generalization of~\cite{mischaikowwanner07a}, Theorem 1.3.
Thus, we begin by recalling some basic definitions and results.

As is indicated in Section~\ref{secintro}, given a continuous
function $u \dvtx  [a,b] \to\mathbb{R}$ and a continuous threshold $\mu
\dvtx [a,b]\to
\mathbb{R}$
we are interested in determining the number of components of the generalized
nodal domain~$N_\mu^\pm$ in terms of a cubical approximation~$Q_{\mu
,M}^\pm$
obtained via sampling at~$M + 1$ points as described in
Definition~\ref{defcubappr}. For suitably chosen discretization
points, and under appropriate regularity and nondegeneracy conditions
on~$u$ one can then expect that the number of components of~$Q_{\mu
,M}^\pm$
and~$N_\mu^\pm$ agree. One only has to be able to verify that the function~$u$
has at most one zero (counting multiplicity) in each of the
intervals~$[x_{k-1},x_k]$, for $k = 1, \ldots, M$. This is accomplished
using the following framework which goes back to Dunnage~\cite{dunnage66a}.

\begin{definition}
\label{defDCO}
A continuous function $u \dvtx  [a,b] \to\mathbb{R}$ has a
\textit{double crossover on the interval $[\alpha,\beta] \subset
[a,b]$}, if
%
\begin{equation} \label{defDCO1}
\sigma\cdot u( \alpha) \ge0   ,\qquad
\sigma\cdot u\biggl( \frac{\alpha+\beta}{2} \biggr) \le0   \quad
\mbox{and}\quad
\sigma\cdot u( \beta) \ge0
\end{equation}
for one choice of the sign~$\sigma\in\{ \pm1 \}$.
\end{definition}

\begin{definition}
\label{defAdmInt}
Let $u \dvtx  [a,b] \to\mathbb{R}$ be a continuous function.
\begin{itemize}
\item The~\textit{dyadic points} in the
interval~$[\alpha,\beta]$ are defined as
\[
d_{n,k} = \alpha+ (\beta- \alpha) \cdot
\frac{k}{2^n}
\qquad  \mbox{for all }
k = 0,\ldots,2^n
 \mbox{ and }
n \in\mathbb{N}_0   .
\]
The~\textit{dyadic subintervals} of $[\alpha,\beta]$ are the intervals~$[d_{n,k},
d_{n,k+1}]$ for all $k = 0,\ldots,2^n - 1$ and $n \in\mathbb{N}_0$.
\item The~interval~$[\alpha,\beta]\subset[a,b]$ is
\textit{admissible} for $u$, if the function~$u$ does not have a double
crossover on any of the dyadic subintervals of~$[\alpha,\beta]$.
\end{itemize}
\end{definition}

It was shown in~\cite{mischaikowwanner07a} that the concept of admissibility
implies the suitability of our nodal domain approximations. More
precisely, the
following is a slight rewording of~\cite{mischaikowwanner07a}, Proposition~2.5.

\begin{proposition}[(Validation criterion)] \label{propVC}
Let $u \dvtx  [a,b] \to\mathbb{R}$ be a continuous function and let $\mu
\dvtx [a,b]\to\mathbb{R}$
be a continuous threshold function. Let~$N_\mu^\pm$
denote the generalized nodal domains of~$u$, and let~$Q_{\mu,M}^\pm$
denote their cubical
approximations as in Definition~\ref{defcubappr}. Furthermore, assume
that the following hold:
\begin{enumerate}[(c)]
\item[(a)] The~function~$u-\mu$ is nonzero at all grid points~$x_k$,
for $k = 0, \ldots, M$.
\item[(b)] The~function~$u-\mu$ has no double zero in~$(a,b)$, that is,
if $x \in(a,b)$ is a zero of~$u$, then~$u-\mu$ attains both
positive and negative function values in every neighborhood
of~$x$.
\item[(c)] For every $k = 1, \ldots, M$, the interval~$[x_{k-1}, x_k]$
between consecutive discretization points is admissible
for~$u-\mu$ in the sense of Definition~\textup{\ref{defAdmInt}}.
\end{enumerate}
Then we have
\[
\beta_0 ( N_\mu^\pm) = \beta_0 ( Q_{\mu,M}^\pm
)   .
\]
\end{proposition}

The~following lemma provides bounds on the probability for admissibility
of a given interval.

\begin{lemma} \label{lemAdmissProb}
Consider a probability space~$(\Omega, \mathcal{F}, \mathbb{P})$,
a continuous threshold function $\mu\dvtx  [a,b] \to\mathbb{R}$, and a
random field
$u \dvtx  [a,b] \times\Omega\to\mathbb{R}$ over~$(\Omega, \mathcal{F},
\mathbb{P})$ such
that~$u(\cdot,\omega)$ is continuous for $\mathbb{P}$-almost all
$\omega\in\Omega$. In addition, assume that~\textup{(A1)}, \textup{(A2)} and~\textup{(A3)} hold.
If $[x,x+\delta] \subset[a,b]$, then
\begin{eqnarray} \label{lemAdmissProb1}
&&\mathbb{P}( [x,x+\delta] \mbox{ is not admissible for } u-\mu
)\nonumber\\[-8pt]\\[-8pt]
&&\qquad  \leq
\frac{4 \mathcal{C}_0(x)}{3} \cdot\delta^3 +
\biggl( \frac{4L}{3} + \frac{8\mathcal{C}_1}{7}\biggr) \cdot\delta
^4,\nonumber
\end{eqnarray}
where $L = \max\{ |\mathcal{C}_0^\prime(y)| \dvtx  y\in[a,b]
\}$.
\end{lemma}

\begin{pf}
If the interval~$I = [x,x+\delta]$ is not admissible, then the
function~$u - \mu$
has a double crossover on one of its dyadic subintervals. If we now
denote the
dyadic points in~$I$ by~$d_{n,k}$ as in Definition~\ref{defAdmInt}, then
together with~(A3) one obtains the estimate
\begin{eqnarray*}
\mathbb{P}\{ I  \mbox{ is not admissible} \} & \le&
\sum_{n=0}^\infty\sum_{k=0}^{2^n - 1} \bigl(
p_{+1} ( d_{n,k}, \delta/ 2^n ) +
p_{-1} ( d_{n,k}, \delta/ 2^n ) \bigr) \\
& \le& \sum_{n=0}^\infty\sum_{k=0}^{2^n - 1} \biggl(
\mathcal{C}_0 ( d_{n,k} ) \cdot\biggl( \frac{\delta}{2^n}
\biggr)^3
+ \mathcal{C}_1 \cdot\biggl( \frac{\delta}{2^n} \biggr)^4 \biggr)
  .
\end{eqnarray*}
Since~$\mathcal{C}_0$ is continuously differentiable, we can define
$L = \max\{ |\mathcal{C}_0^\prime(y)| \dvtx  y \in[a,b] \}$, and the definition
of the dyadic points implies
\[
\mathcal{C}_0 ( d_{n,k} ) \le\mathcal{C}_0(x) + L \cdot
(d_{n,k} - x)
\le
\mathcal{C}_0(x) + L \delta  .
\]
This finally furnishes
\begin{eqnarray*}
\mathbb{P}\{ I  \mbox{ is not admissible} \} & \le&
\sum_{n=0}^\infty\sum_{k=0}^{2^n - 1} \biggl(
\bigl( \mathcal{C}_0(x) + L \delta\bigr) \cdot
\biggl( \frac{\delta}{2^n} \biggr)^3
+ \mathcal{C}_1 \cdot\biggl( \frac{\delta}{2^n} \biggr)^4 \biggr)
\\
& = & \frac{4 \mathcal{C}_0(x)}{3} \cdot\delta^3 +
\biggl( \frac{4 L}{3} + \frac{8 \mathcal{C}_1}{7} \biggr) \cdot
\delta^4
  .
\end{eqnarray*}
\upqed\end{pf}
Combining Proposition~\ref{propVC}, Lemma~\ref{lemAdmissProb},
and restricting to the leading order term in~(\ref{lemAdmissProb1}),
one obtains
%
\begin{equation} \label{heur1}
\mathbb{P}\{ \beta_0 ( N_\mu^\pm) = \beta_0
( Q_{\mu,M}^\pm) \}   \ge
1 -   \frac{4}{3} \cdot
\sum_{k=1}^M \mathcal{C}_0(x_{k-1}) \cdot( x_k - x_{k-1}
)^3   .
\end{equation}
Clearly, the resulting bound depends on the location of the
sampling points, which suggests maximizing the bound to optimize the location.

We first provide a heuristic argument for this optimal location, and present
the precise result afterwards. One can show that for arbitrary
nonnegative numbers~$\delta_1, \ldots, \delta_M \ge0$ the inequality
\[
\sum_{k=1}^M \delta_k^3   \ge
\frac{1}{M^2} \cdot\Biggl( \sum_{k=1}^M \delta_k \Biggr)^3
\]
holds, with equality if and only if $\delta_1 = \delta_2 = \cdots=
\delta_M$.
Applying this inequality to the sum in the right-hand side of~(\ref{heur1}),
implies
%
\begin{equation} \label{heur2}
\quad \sum_{k=1}^M \mathcal{C}_0(x_{k-1}) \cdot( x_k - x_{k-1}
)^3
\ge
\frac{1}{M^2} \cdot\Biggl( \sum_{k=1}^M \sqrt[3]{\mathcal{C}_0(x_{k-1})}
\cdot
( x_k - x_{k-1} ) \Biggr)^3
\end{equation}
with equality if and only if
%
\begin{eqnarray} \label{heur3}
\sqrt[3]{\mathcal{C}_0(x_{k-1})} \cdot( x_k - x_{k-1} ) =
\sqrt[3]{\mathcal{C}_0(x_{\ell-1})} \cdot( x_\ell- x_{\ell-1}
)\nonumber\\[-8pt]\\[-8pt]
 \eqntext{\mbox{for all } k, \ell= 1, \ldots, M   .}
\end{eqnarray}
For large~$M$, the sum on the right-hand side of~(\ref{heur2})
converges to the integral of~$\mathcal{C}_0^{1/3}$ over~$[a,b]$.
The~motivation for Theorem~\ref{thmTopSamp} is now clear:
Condition~(\ref{heur3}) suggests that for $M \to\infty$,
the optimal estimate can be achieved by choosing the sampling points in an
equi-$\mathcal{C}_0^{1/3}$-area fashion, since the term $\mathcal
{C}_0(x_{k-1})^{1/3}
(x_k - x_{k-1})$ approximates the intergral of~$\mathcal{C}_0^{1/3}$
over~$[x_{k-1}, x_k]$. This heuristic forms the basis for the following
proof of our first main result.

\begin{pf*}{Proof of Theorem~\ref{thmTopSamp}}
Let $\delta_{\max} := \max_{k=1,\ldots,M} |x_k - x_{k-1}|$, and define
the positive number $m := \min_{x \in[a,b]} \mathcal{C}_0(x)^{1/3} >
0$. Furthermore,
let $L := \max_{x \in[a,b]} |d\mathcal{C}_0^{1/3} / dx|$. Then the mean
value theorem
readily furnishes
%
\begin{equation} \label{thmTopSamp2}
\bigg| \sqrt[3]{\mathcal{C}_0(x_{k-1})} \cdot( x_k - x_{k-1}
) -
\int_{x_{k-1}}^{x_k} \sqrt[3]{\mathcal{C}_0(x)}  \, dx \bigg|
\le
L ( x_k - x_{k-1} )^2
\end{equation}
for all $k=1,\ldots,M$. Due to the choice of the sampling points we further
have
%
\begin{equation} \label{thmTopSamp3}
m \cdot(x_k - x_{k-1})   \le
\int_{x_{k-1}}^{x_k} \sqrt[3]{\mathcal{C}_0(x)}  \, dx   =
\frac{1}{M} \cdot\underbrace{\int_a^b \sqrt[3]{\mathcal{C}_0(x)}
\,dx}_{=: K}   ,
\end{equation}
which in turn implies
%
\begin{equation} \label{thmTopSamp4}
0   <   x_k - x_{k-1}   \le
\delta_{\max}   \le
\frac{K}{m \cdot M}
\qquad  \mbox{for all }
k=1,\ldots,M   .
\end{equation}
Applying Lemma~\ref{lemAdmissProb} to every subinterval formed by adjacent
sampling points, we now obtain together with~(\ref{thmTopSamp2}),
(\ref{thmTopSamp3}) and~(\ref{thmTopSamp4}) the estimate
\begin{eqnarray*}
&&1 - \mathbb{P}\{ \beta_0(N_\mu^\pm) =\beta_0(Q_{\mu,M}^\pm)
\}\\
&&\qquad  \le \frac{4}{3} \sum_{k=1}^M \mathcal{C}_0(x_{k-1}) \cdot(x_k -
x_{k-1})^3
+ \mathcal{C}_2 \sum_{k=1}^M (x_k - x_{k-1})^4 \\
&&\qquad  \le \frac{4}{3} \cdot\sum_{k=1}^M \biggl( \frac{K}{M} + L
( x_k - x_{k-1})^2 \biggr)^3
+ \frac{\mathcal{C}_2 K^4}{m^4 M^3} \\
&&\qquad  \le \frac{4}{3} \cdot\sum_{k=1}^M \biggl( \frac{K}{M} +
\frac{L K^2}{m^2 M^2} \biggr)^3
+ \frac{\mathcal{C}_2 K^4}{m^4 M^3} \\
&&\qquad  =  \frac{4 K^3}{3 M^2 } + O\biggl( \frac{1}{M^3} \biggr)
\end{eqnarray*}
for some constant~$\mathcal{C}_2 \ge0$. This is exactly~(\ref{thmTopSamp1}).
\end{pf*}

\section{Asymptotics of sign-change probabilities} \label{secPTool}
Theorem~\ref{thmCorrSamp} can be viewed as a special case of
Theorem~\ref{thmTopSamp}. The~content lies in the fact that under the
assumption of a Gaussian random field, the function $\mathcal{C}_0$
can be
explicitly computed. However, this requires a quantitative understanding
of the asymptotic behavior of sign-distribution probabilities of
parameter-dependent Gaussian random variables, which is the focus of this
section.

More precisely, let~$T(\delta) = (T_1(\delta),\ldots,T_n(\delta))^t
\in\mathbb{R}^n$
denote a one-parameter family of $\mathbb{R}^n$-valued random Gaussian
variables over
a probability space\break  $(\Omega,\mathcal{F},\mathbb{P})$, indexed by
$\delta> 0$, and
choose a sign
sequence $(s_1,\ldots,s_n) \in\{ \pm1 \}^n$. Furthermore, let~$\tau
(\delta)
\in\mathbb{R}^3$ denote an arbitrary threshold vector. We are
interested in the
precise asymptotic behavior as $\delta\to0$ of the probability
%
\begin{equation} \label{PTool1}
P(\delta) = \mathbb{P}\bigl\{ s_j \bigl( T_j(\delta) - \tau
_j(\delta)
\bigr) \ge0
 \mbox{ for all }  j = 1,\ldots,n \bigr\}   .
\end{equation}
The~following result is an extension
of~(\cite{mischaikowwanner07a}, Proposition~4.1)
which dealt only with the special case $\tau\equiv0$.

\begin{proposition} \label{propPToolA}
Let $(s_1,\ldots,s_n) \in\{ \pm1 \}^n$ denote a fixed sign sequence,
and consider one-parameter families of a threshold vector $\tau(\delta)
\in\mathbb{R}^3$ and an $\mathbb{R}^n$-valued random Gaussian
variable~$T(\delta)$
over a probability space~$(\Omega,\mathcal{F},\mathbb{P})$, for
$\delta> 0$.
Assume that the following hold:
\begin{enumerate}[(iii)]
\item[(i)] For each $\delta> 0$, assume that the Gaussian random
variable~$T(\delta)$ has mean~$0 \in\mathbb{R}^n$ and a positive definite
covariance matrix~$C(\delta) \in\mathbb{R}^{n \times n}$, whose positive
eigenvalues are given by $0 < \lambda_1(\delta) \le\cdots\le
\lambda_n(\delta)$. The~corresponding orthonormalized eigenvectors
are denoted by $v_1(\delta), \ldots, v_n(\delta)$.
\item[(ii)] There exists a vector~$\bar{v}_1 = (\bar{v}_{11}, \ldots,
\bar{v}_{1n})^t \in\mathbb{R}^n$ such that $v_1(\delta) \to\bar{v}_1$
as \mbox{$\delta\to0$}, and $s_j \cdot\bar{v}_{1j} > 0$ for all $j =
1,\ldots,n$.
\item[(iii)] The~quotient $\lambda_1(\delta) / \lambda_k(\delta)$ converges
to~$0$ as $\delta\to0$, for all $k = 2,\ldots,n$.
\item[(iv)] There exists a vector $\alpha= (\alpha_1, \ldots,
\alpha_n)^t
\in\mathbb{R}^n$ such that
%
\begin{equation} \label{propPToolA1}
\lim_{\delta\to0}
\frac{\tau(\delta) \cdot v_k(\delta)}{\lambda_k(\delta)^{1/2}}
= \alpha_k
\qquad  \mbox{for all }
k = 1,\ldots,n   .
\end{equation}
\end{enumerate}
Furthermore, for~$\alpha$ as above define
%
\begin{equation} \label{propPToolA2}
S_\alpha= \frac{2}{2^{n/2} \cdot\Gamma(n/2)} \cdot
e^{-\sum_{k=2}^n \alpha_k^2 / 2} \cdot\int_{\alpha_1}^\infty
( s - \alpha_1 )^{n-1} e^{-s^2 / 2} \,  ds   .
\end{equation}
Then the probability~$P(\delta)$ defined in~\textup{(\ref{PTool1})} satisfies
%
\begin{equation} \label{propPToolA3}
\lim_{\delta\to0} P(\delta) \cdot
\sqrt{\frac{\det C(\delta)}{\lambda_1(\delta)^n}}   =
\frac{\Gamma(n/2) \cdot S_\alpha}{2 \cdot\pi^{n/2} \cdot(n-1)!}
\cdot\Bigg| \prod_{j=1}^n \bar{v}_{1j} \Bigg|^{-1}   .
\end{equation}
\end{proposition}

For specific values of~$n$, the integral in~(\ref{propPToolA1}) can
be simplified further. For our one-dimensional application, we need
the case $n = 3$, which is the subject of the following remark.

\begin{remark} \label{remPToolA}
Recall that $\Gamma(1/2) = \pi^{1/2}$, $\Gamma(1) = 1$, and
$\Gamma(t+1) = t \Gamma(t)$ for $t > 0$. Furthermore, notice that
$S_\alpha= 1$ for $\alpha= 0 \in\mathbb{R}^n$. In addition,
for $n = 3$ one can readily verify that
\begin{eqnarray} \label{propPToolA4}
S_\alpha&=& \frac{2^{1/2}}{\pi^{1/2}} \cdot
e^{-( \alpha_2^2 + \alpha_3^2 ) / 2} \nonumber\\[-8pt]\\[-8pt]
&&{}\times \biggl(
-\alpha_1 e^{-\alpha_1^2 / 2} + ( 1 + \alpha_1^2 )
\cdot\int_{\alpha_1}^\infty e^{-s^2 / 2}  \, ds \biggr)   .\nonumber
\end{eqnarray}
\end{remark}

\begin{pf}
Define the diagonal matrix~$S = (s_{i} \delta_{ij})_{i,j=1,\ldots,n}$,
where~$\delta_{ij}$ denotes the Kronecker delta, and let~$Z_+ = \{ z
\in\mathbb{R}^n \dvtx  z_j \ge0 \mbox{ for } j=1,\ldots,n \}$. Finally, let
\[
D(\delta) = \lambda_1(\delta) \cdot S C(\delta)^{-1} S
\]
and
\[
d(\delta) = \frac{1}{\lambda_1(\delta)^{1/2}} \cdot S \tau(\delta
)   .
\]
Using the density of the Gaussian distribution of~$T(\delta)$
according to~(\cite{bauer96a}, Theorem~30.4), which exists
since~$C(\delta)$ is positive definite, in combination with a simple
rescaling and shifting of the coordinate system, the probability
in~(\ref{PTool1}) can be rewritten as
\begin{eqnarray*}
P(\delta) & = & \frac{(2\pi)^{-n/2}}{\sqrt{\det C(\delta)}} \cdot
\int_{S\tau(\delta) + Z_+} e^{-z^t S C(\delta)^{-1} S z/2}   \,dz \\
& = & \sqrt{\frac{\lambda_1(\delta)^n}{2^n \pi^n \det C(\delta)}}
\cdot
\int_{Z_+} e^{-(z + d(\delta) )^t D(\delta)
(z + d(\delta) ) / 2}   \,dz   .
\end{eqnarray*}
According to our assumptions, the eigenvalues~$\mu_1(\delta), \ldots,
\mu_n(\delta)$ of the matrix~$D(\delta)$ are given by
\[
\mu_1(\delta) = 1   \quad
 \mbox{and}\quad
\mu_k(\delta) = \frac{\lambda_1(\delta)}{\lambda_k(\delta)}
\qquad  \mbox{for }
k = 2,\ldots,n   ,
\]
with corresponding orthonormalized eigenvectors~$w_k(\delta) =
Sv_k(\delta)$, for $k = 1,\break \ldots,n$. Now let~$B(\delta)$ denote the
orthogonal matrix with columns~$w_1(\delta), \ldots,\break w_n(\delta)$ and
introduce the change of variables $z = B(\delta)\zeta$. Moreover,
let
\[
Z(\zeta_1,\delta) = \Biggl\{ (\zeta_2,\ldots,\zeta_n)   \dvtx
\sum_{k=1}^n \zeta_k w_k(\delta) \in Z_+ \Biggr\}
\subset\mathbb{R}^{n-1}
\]
define real numbers~$\eta_1(\delta), \ldots, \eta_n(\delta)$ by
\[
\eta_k(\delta) = S\tau(\delta) \cdot w_k(\delta) =
\tau(\delta) \cdot v_k(\delta)
\qquad  \mbox{for }
k = 1,\ldots,n   ,
\]
and let
\[
I(\zeta_1,\delta) =
\int_{Z(\zeta_1,\delta)} \exp\Biggl( -\sum_{k=1}^n
\frac{\mu_k(\delta)}{2}
\biggl( \zeta_k + \frac{\eta_k(\delta)}{\lambda_1(\delta)^{1/2}}
\biggr)^2 \Biggr) \,  d(\zeta_2,\ldots,\zeta_n)   .
\]
Due to~(ii) and the definition of the signs~$s_k$, the
eigenvector~$w_1(\delta)$ has strictly positive components for all
sufficiently small $\delta> 0$, and therefore the identity
\[
\bigl( z + d(\delta) \bigr)^t D(\delta) \bigl( z + d(\delta)
\bigr) = \sum_{k=1}^n \mu_k(\delta)
\biggl( \zeta_k + \frac{\eta_k(\delta)}{\lambda_1(\delta)^{1/2}}
\biggr)^2
\]
implies
%
\begin{eqnarray}\label{propPToolA5}
&&\int_{Z_+} e^{-( z + d(\delta) )^t D(\delta)
( z + d(\delta) ) / 2}  \, dz\nonumber\\
 &&\quad  =
\int_{B(\delta)^{-1}Z_+} \exp\Biggl( -\sum_{k=1}^n
\frac{\mu_k(\delta)}{2}
\biggl( \zeta_k + \frac{\eta_k(\delta)}{\lambda_1(\delta)^{1/2}}
\biggr)^2 \Biggr) \,d\zeta\\
&&\quad  =  \int_0^\infty I(\zeta_1,\delta)\,   d\zeta_1   .\nonumber
\end{eqnarray}
From the definition of~$I(\zeta_1,\delta)$, one can easily deduce
\[
I(\zeta_1,\delta) = \zeta_1^{n-1} \cdot
\int_{Z(1,\delta)} \exp\Biggl( -\sum_{k=1}^n \frac{\mu_k(\delta)}{2}
\biggl( \zeta_1 \xi_k + \frac{\eta_k(\delta)}{\lambda_1(\delta)^{1/2}}
\biggr)^2 \Biggr) \,  d(\xi_2,\ldots,\xi_n)   ,
\]
where we define~$\xi_1 = 1$. This representation furnishes
for all $\zeta_1 > 0$ and $\delta> 0$ the estimate
%
\begin{equation} \label{propPToolA6}
I(\zeta_1,\delta) \le\zeta_1^{n-1} \cdot
\operatorname{vol}_{n-1} (Z(1,\delta)) \cdot
e^{-( \zeta_1 + \eta_1(\delta) \lambda_1(\delta)^{-1/2}
)^2 / 2}   .
\end{equation}
Again according to (ii), the $(n-1)$-dimensional volume of the
simplex~$Z(1,\delta)$ converges to the $(n-1)$-dimensional volume of
the simplex
\[
\widetilde{Z} = \{ z \in Z_+   \dvtx
( z - S\bar{v}_1 , S\bar{v}_1 ) = 0 \}
\subset\mathbb{R}^n   ,
\]
which can be computed as
\[
\operatorname{vol}_{n-1} ( \widetilde{Z} ) =
\frac{1}{(n-1)!} \cdot
\Bigg| \prod_{j=1}^n \bar{v}_{1j} \Bigg|^{-1}   .
\]
Now let~$\zeta_1 > 0$ be arbitrary, but fixed. Notice that since we
did not make any assumptions about the asymptotic behavior of the
eigenvectors $w_2(\delta), \ldots, w_n(\delta)$ for $\delta\to0$,
the sets~$Z(1,\delta)$ do not have to converge. Yet, (ii)
yields the existence of a compact subset~$K \subset\mathbb{R}^{n-1}$ such
that $Z(1,\delta) \subset K$ for all sufficiently small $\delta>
0$. Furthermore, we have
\begin{eqnarray*}
\sum_{k=1}^n \mu_k(\delta) \biggl( \zeta_1 \xi_k +
\frac{\eta_k(\delta)}{\lambda_1(\delta)^{1/2}} \biggr)^2 & = &
\zeta_1^2 + \frac{2 \zeta_1 \eta_1(\delta)}{\lambda_1(\delta)^{1/2}}
+ \frac{\eta_1(\delta)^2}{\lambda_1(\delta)} + \sum_{k=2}^n
\frac{\zeta_1^2 \xi_k^2 \lambda_1(\delta)}{\lambda_k(\delta)} \\
& &   +   2 \sum_{k=2}^n \frac{\zeta_1 \xi_k \eta_k(\delta)
\lambda_1(\delta)^{1/2}}{\lambda_k(\delta)}
+ \sum_{k=2}^n \frac{\eta_k(\delta)^2}{\lambda_k(\delta)} \\
& \to& \zeta_1^2 + 2 \zeta_1 \alpha_1 + \alpha_1^2
+ \sum_{k=2}^n \alpha_k^2
\end{eqnarray*}
as $\delta\to0$. Due to (iii) and (iv), this convergence
is uniform on~$K$. Therefore, we have
%
\begin{eqnarray}
\lim_{\delta\to0} I(\zeta_1,\delta) =
\zeta_1^{n-1} \cdot
\operatorname{vol}_{n-1} ( \widetilde{Z} ) \cdot
e^{-( \zeta_1 + \alpha_1 )^2 / 2} \cdot
e^{-( \alpha_2^2 + \cdots+ \alpha_n^2 ) / 2}\nonumber\\
 \eqntext{\mbox{for all }\zeta_1 > 0   .}
\end{eqnarray}
Due to~(\ref{propPToolA6}) and $\operatorname{vol}_{n-1}(Z(1,\delta)) \to
\operatorname{vol}_{n-1}(\widetilde{Z})$, we can now apply the dominated
convergence theorem to pass to the limit $\delta\to0$
in~(\ref{propPToolA5}), and this furnishes
\begin{eqnarray*}
& & \lim_{\delta\to0} \int_{Z_+} e^{-( z + d(\delta) )^t
D(\delta) ( z + d(\delta) ) / 2}  \, dz \\
&&\qquad  =  \operatorname{vol}_{n-1} ( \widetilde{Z} ) \cdot
e^{-( \alpha_2^2 + \cdots+ \alpha_n^2 ) / 2} \cdot
\int_0^\infty\zeta_1^{n-1} e^{-( \zeta_1 + \alpha_1
)^2 / 2}
 \, d\zeta_1 \\
& &\qquad =  \operatorname{vol}_{n-1} ( \widetilde{Z} ) \cdot
e^{-( \alpha_2^2 + \cdots+ \alpha_n^2 ) / 2} \cdot
\int_{\alpha_1}^\infty( s - \alpha_1 )^{n-1} e^{-s^2 / 2}
 \, ds.
\end{eqnarray*}
\upqed\end{pf}

We close this section with a corollary to Proposition~\ref{propPToolA}.
In our applications of the above result, we are not only interested in the
asymptotic behavior of~$P(\delta)$ as defined in~(\ref{PTool1}), that is,
for the fixed sign sequence~$(s_1, \ldots, s_n)$, but also in the
corresponding probability for the negative sign sequence~$(-s_1, \ldots
, -s_n)$.

More precisely, if~$T(\delta) = (T_1(\delta),\ldots,T_n(\delta))^t
\in\mathbb{R}^n$
denotes again a one-parameter family of $\mathbb{R}^n$-valued random Gaussian
variables
over a probability space~$(\Omega,\mathcal{F},\mathbb{P})$, indexed
by $\delta>
0$, and if we
choose both a sign sequence $(s_1,\ldots, s_n) \in\{ \pm1 \}^n$ and a
one-parameter family~$\tau(\delta) \in\mathbb{R}^n$ of threshold vectors,
then we are interested in the asymptotic behavior as $\delta\to0$ of the
probability
\begin{eqnarray}\label{PTool2}
P^{\pm}(\delta) & = & \mathbb{P}\bigl\{ s_j \bigl( T_j(\delta) -
\tau_j(\delta)\bigr) \ge0 \mbox{ for all }  j = 1,\ldots,n \bigr\} \nonumber\\[-8pt]\\[-8pt]
& & {} +    \mathbb{P}\bigl\{ s_j \bigl( T_j(\delta) - \tau
_j(\delta)
\bigr) \le0
 \mbox{ for all }  j = 1,\ldots,n \bigr\}   .\nonumber
\end{eqnarray}
This is the subject of the following corollary.

\begin{corollary} \label{corPToolA}
Let $(s_1,\ldots,s_n) \in\{ \pm1 \}^n$ denote a fixed sign sequence,
let $\tau(\delta) \in\mathbb{R}^n$ denote a threshold vector,
and consider a one-parameter family~$T(\delta)$, $\delta> 0$, of
$\mathbb{R}^n$-valued random Gaussian variables over a probability
space~$(\Omega,\mathcal{F},\mathbb{P})$ which satisfies all the
assumptions of
Proposition~\textup{\ref{propPToolA}}. Then the probability~$P^{\pm}(\delta)$
defined in~\textup{(\ref{PTool2})} satisfies
%
\begin{equation} \label{corPToolA1}
\lim_{\delta\to0} P^{\pm}(\delta) \cdot
\sqrt{\frac{\det C(\delta)}{\lambda_1(\delta)^n}}   =
\frac{\Gamma(n/2) \cdot S_\alpha^{\pm}}{2 \cdot\pi^{n/2}
\cdot(n-1)!} \cdot\Bigg| \prod_{j=1}^n \bar{v}_{1j} \Bigg|^{-1}
  ,
\end{equation}
where~$S_\alpha^{\pm} = S_\alpha+ S_{-\alpha}$, with~$\alpha$ as
in~\textup{(\ref{propPToolA1})} and~$S_\alpha$ as in~\textup{(\ref{propPToolA2})}.
Moreover, for the special case $n = 3$ one obtains
%
\begin{equation} \label{corPToolA2}
S_\alpha^{\pm} =
2   e^{-( \alpha_2^2 + \alpha_3^2 ) / 2} \cdot
( 1 + \alpha_1^2 )   .
\end{equation}
\end{corollary}

\begin{pf}
One only has to apply Proposition~\ref{propPToolA} twice---first
with the given
sign vector~$(s_1, \ldots, s_n)$, and then with the sign
vector~$(-s_1, \ldots, -s_n)$. Notice that in the latter case, we have
to use the eigenvector~$-v_1(\delta)$ instead of~$v_1(\delta)$, which leads
to~$-\alpha_k$ instead of~$\alpha_k$ in~(\ref{propPToolA1});
everything else
remains unchanged. This immediately implies~(\ref{corPToolA1}). As
for~(\ref{corPToolA2}), one only has to notice that
\[
\int_{\alpha_1}^{\infty} e^{-s^2 / 2}  \, ds +
\int_{-\alpha_1}^{\infty} e^{-s^2 / 2} \,  ds =
\int_{-\infty}^{\infty} e^{-s^2 / 2}  \, ds =
\sqrt{2\pi}
\]
and employ Remark~\ref{remPToolA}.
\end{pf}

\section{Sampling based on spatial correlations} \label{secspatial}
The~goal of this section is the proof of Theorem~\ref{thmCorrSamp}. To do
this, we need to relate the spatial correlation function~$R$ to local
probability asymptotics. For this, we use the following lemma.

\begin{lemma}
\label{lemCorrSamp}
Consider a Gaussian random field~$u \dvtx  [a,b] \times\Omega\to\mathbb
{R}$ satisfying
\textup{(G1)} and \textup{(G2)}.
For $x \in[a,b)$ and sufficiently small values of $\delta> 0$, define the
random vector~$T(\delta) = (T_1(\delta), T_2(\delta), T_3(\delta
))^t$ via
%
\begin{equation} \label{lemCorrSamp1}
T_1(\delta) = u( x )   ,\qquad
T_2(\delta) = u\biggl( x + \frac{\delta}{2} \biggr) \quad
\mbox{and}\quad
T_3(\delta) = u( x + \delta)   .
\end{equation}
Then~$T$ is a centered Gaussian random variable with positive definite
covariance matrix~$C(\delta)$. Moreover, if we denote the eigenvalues
of~$C(\delta)$ by $0 < \lambda_1(\delta) \le\lambda_2(\delta) \le
\lambda_3(\delta)$, then
\begin{eqnarray*}
\lambda_1(\delta) & = & \frac{\det\mathcal{R}(x)}{96   \mathcal
{R}_{3,3}^m(x)}
\cdot\delta^4 + O( \delta^5 ), \\
\lambda_2(\delta) & = & \frac{\mathcal{R}_{3,3}^m(x)}{2   R_{0,0}(x)}
\cdot\delta^2 + O( \delta^3 ), \\
\lambda_3(\delta) & = & 3   R_{0,0}(x) + O( \delta),
\end{eqnarray*}
where we use the notation introduced in~\textup{(\ref{assRFcorr1})}, \textup{(\ref{assRFcorr2})}
and~\textup{(\ref{defRFcorr})}.
In addition, we can choose the normalized eigenvectors~$v_1(\delta)$,
$v_2(\delta)$ and~$v_3(\delta)$ corresponding to these eigenvalues in
such a way that
\begin{eqnarray*}
\lim_{\delta\to0} v_1(\delta) &=& \frac{1}{\sqrt{6}}
\pmatrix{  1 \cr -2 \cr 1}   ,\qquad
\lim_{\delta\to0} v_2(\delta) = \frac{1}{\sqrt{2}}
\pmatrix{ 1 \cr 0 \cr -1
} ,\\
\lim_{\delta\to0} v_3(\delta) &=& \frac{1}{\sqrt{3}}
\pmatrix{  1 \cr 1 \cr 1}.
\end{eqnarray*}
Finally, for a $C^3$-function $\mu\dvtx  [a,b] \to\mathbb{R}$ define the vector
$\tau(\delta) = (\tau_1(\delta),\tau_2(\delta),\break \tau_3(\delta))^t$
via
%
\begin{equation} \label{lemCorrSamp2}
\tau_1(\delta) = \mu( x )   ,\qquad
\tau_2(\delta) = \mu\biggl( x + \frac{\delta}{2} \biggr)
\quad
\mbox{and}\quad
\tau_3(\delta) = \mu( x + \delta)   .
\end{equation}
Then
\begin{eqnarray*}
\tau(\delta) \cdot v_1(\delta) & = &
\frac{\mathcal{R}_{3,1}^m(x) \mu(x) - \mathcal{R}_{3,2}^m(x) \mu
^\prime(x) +
\mathcal{R}_{3,3}^m(x)\mu^{\prime\prime}(x)}{4\sqrt{6}
\mathcal{R}_{3,3}^m(x)}
\cdot\delta^2 + O( \delta^3 )   , \\
\tau(\delta) \cdot v_2(\delta) & = &
\frac{R_{1,0}(x) \mu(x) - R_{0,0}(x) \mu^\prime(x)}{\sqrt{2}
R_{0,0}(x)} \cdot\delta+ O( \delta^2 )   , \\
\tau(\delta) \cdot v_3(\delta) & = &
\sqrt{3} \cdot\mu(x) + O( \delta)   .
\end{eqnarray*}
\end{lemma}

\begin{pf}
Due to our assumptions on~$u$, the vector~$T(\delta)$ is normally
distributed with mean~$0 \in\mathbb{R}^3$ and covariance
matrix~$C(\delta)
\in\mathbb{R}^{3 \times3}$ given by
\[
C(\delta) = \pmatrix{
r(0,0) & r(0,\delta/2) & r(0,\delta) \vspace*{2pt}\cr
r(0,\delta/2) & r(\delta/2,\delta/2) & r(\delta/2,\delta) \vspace*{2pt}\cr
r(0,\delta) & r(\delta/2,\delta) & r(\delta,\delta)
} ,
\]
where we use the abbreviation
\[
r(\delta_1,\delta_2)  =  R(x+\delta_1, x+\delta_2)   .
\]
For $(\delta_1,\delta_2) \to0$, the function~$r$ can be expanded as
\begin{eqnarray*}
r(\delta_1,\delta_2) & = & R_{0,0}(x) + R_{1,0}(x) \delta_1 +
R_{1,0}(x) \delta_2 + \frac{R_{2,0}(x)}{2} \delta_1^2 +
R_{1,1}(x) \delta_1 \delta_2\\
&&{} + \frac{R_{2,0}(x)}{2} \delta_2^2
 + \frac{R_{3,0}(x)}{6} \delta_1^3 +
\frac{R_{2,1}(x)}{2} \delta_1^2 \delta_2 +
\frac{R_{2,1}(x)}{2} \delta_1 \delta_2^2\\
&&{} +
\frac{R_{3,0}(x)}{6} \delta_2^3 +
O( |(\delta_1,\delta_2)|^4 )   ,
\end{eqnarray*}
where the~$R_{k,\ell}$ where defined in~(\ref{assRFcorr1}). Furthermore,
(G2) implies that we have the
strict inequalities
\[
R_{0,0}(x) > 0   ,\qquad
\mathcal{R}_{3,3}^m(x) > 0   \quad   \mbox{as well as}\quad
\det\mathcal{R}(x) > 0   .
\]
These strict inequalities ensure that in all of the
expansions derived below the leading order coefficients are positive.

Using the above expansion of~$r$, the determinant of the covariance
matrix~$C(\delta)$ of the random vector~$T(\delta)$ can be written as
\[
\det C(\delta)   =   \tfrac{1}{64} \cdot\det\mathcal{R}(x)
\cdot\delta^6 + O( \delta^7 )  ,
\]
that is, the covariance matrix is positive definite for sufficiently
small $\delta> 0$. Furthermore, by applying the
Newton polygon method~\cite{sidorovetal02a,vainbergtrenogin74a}
to the
characteristic polynomial~$\det( C(\delta) - \lambda I)$ it can be
shown that in the limit $\delta\to0$ the three
eigenvalues~$\lambda_k(\delta)$, for $k=1,2,3$, of~$C(\delta)$ are
given by
the expansions in the formulation of Lemma~\ref{lemCorrSamp}.

We now turn our attention to the asymptotic statements concerning the
eigenvectors of the covariance matrix. According to the form
of $C(\delta)$,
we have
\[
\lim_{\delta\to0} C(\delta) =
\pmatrix{
R_{0,0}(x) & R_{0,0}(x) & R_{0,0}(x)\vspace*{2pt} \cr
R_{0,0}(x) & R_{0,0}(x) & R_{0,0}(x)\vspace*{2pt} \cr
R_{0,0}(x) & R_{0,0}(x) & R_{0,0}(x)
},
\]
where the limit has a double eigenvalue~$0$, as well as
the simple eigenvalue $3 R_{0,0}(x)$ with normalized
eigenvector $(1,1,1)^t / 3^{1/2}$. Due to standard results on the
perturbation of simple eigenvalues and corresponding
eigenvectors~\cite{wilkinson88a}, this implies that~$v_3(\delta)$ can
be chosen as in the formulation of the lemma.

In order to determine the asymptotic behavior of the eigenvector
corresponding to~$\lambda_1$, we consider the adjoint of the covariance
matrix, whose expansion is given by
\[
\operatorname{adj}   C(\delta) = \frac{\mathcal{R}_{3,3}^m(x)}{4} \cdot
\pmatrix{
1 & -2 & 1 \vspace*{2pt}\cr
-2 & 4 & -2\vspace*{2pt} \cr
1 & -2 & 1
} \cdot\delta^2
+ O( \delta^3 )  .
\]
The~constant coefficient matrix has the double eigenvalue~$0$, as well
as the positive eigenvalue~$6$ with associated unnormalized
eigenvector~$(1,-2,1)^t$. Since the eigenspace for the largest
eigenvalue of the adjoint matrix coincides with the eigenspace for the
eigenvalue~$\lambda_1(\delta)$ of~$C(\delta)$, the simplicity of these
eigenvalues shows that we can choose a normalized
eigenvector~$v_1(\delta)$
for~$\lambda_1(\delta)$ with $v_1(\delta) \to(1, -2, 1)^t
/ 6^{1/2}$ for $\delta\to0$. Finally, the orthogonality of the
three eigenvectors shows that we can choose a normalized
eigenvector~$v_2(\delta)$ for~$\lambda_2(\delta)$ with $v_2(\delta)
\to
(1, 0, -1)^t / 2^{1/2}$ for $\delta\to0$.

We now turn our attention to the asymptotics of the inner
products~$\tau(\delta) \cdot v_k(\delta)$. Since~$\mu$ is a
$C^3$-function, we can write
\[
\tau(\delta)   =
\mu(x)
\pmatrix{  1 \cr 1 \cr 1
} +
\frac{\mu^\prime(x)}{2}
\pmatrix{  0 \cr 1 \cr 2
}
\cdot\delta+
\frac{\mu^{\prime\prime}(x)}{8}
\pmatrix{  0 \cr 1 \cr 4
}
\cdot\delta^2 + O( \delta^3 )   ,
\]
and this representation immediately furnishes
$\tau(\delta) \cdot v_3(\delta) \to3^{1/2} \cdot\mu(x)$ as
$\delta\to0$. The~statements concerning $\tau(\delta) \cdot
v_1(\delta)$ and $\tau(\delta) \cdot v_2(\delta)$ are more involved,
and rely on expansions of the eigenvectors in terms of~$\delta$.

As for the first eigenvector, write
$v_1(\delta) = (v_{1,1}(\delta), v_{1,2}(\delta), v_{1,3}(\delta))^t$,
and consider the functions
\[
w_{1,k}(\delta) = -\frac{2}{\sqrt{6} \cdot v_{1,2}(\delta)} \cdot
v_{1,k}(\delta)
\qquad  \mbox{for }
k = 1,2,3   .
\]
Then the vector $w_1(\delta) = (w_{1,1}(\delta), w_{1,2}(\delta),
w_{1,3}(\delta))^t$ is defined for sufficiently small $\delta> 0$,
and for these~$\delta$ we have
\[
w_{1,2}(\delta) = -\tfrac{2}{\sqrt{6}}
\quad  \mbox{as well as}\quad
\bigl( C(\delta) - \lambda_1(\delta) I \bigr) w_1(\delta) = 0   .
\]
Using the abbreviation $C(\delta) = (c_{i,j}(\delta))_{i,j=1,2,3}$, the
latter system is equivalent to
\begin{eqnarray*}
\bigl( c_{1,1}(\delta) - \lambda_1(\delta) \bigr) w_{1,1}(\delta
) +
c_{1,3}(\delta) w_{1,3}(\delta) & = & \tfrac{2}{\sqrt{6}}
\cdot c_{1,2}(\delta)   , \\
c_{3,1}(\delta) w_{1,1}(\delta) + \bigl( c_{3,3}(\delta) -
\lambda_1(\delta) \bigr) w_{1,3}(\delta) & = & \tfrac{2}{\sqrt{6}}
\cdot c_{3,2}(\delta)   ,
\end{eqnarray*}
which immediately implies
\begin{eqnarray*}
w_{1,1}(\delta) & = &
\frac{2}{\sqrt{6}} \cdot\frac{( c_{3,3}(\delta) -
\lambda_1(\delta) ) c_{1,2}(\delta) - c_{3,2}(\delta)
c_{1,3}(\delta)}{( c_{1,1}(\delta) -
\lambda_1(\delta) ) ( c_{3,3}(\delta) -
\lambda_1(\delta) ) - c_{1,3}(\delta) c_{3,1}(\delta)}
  , \\
w_{1,3}(\delta) & = &
\frac{2}{\sqrt{6}} \cdot\frac{( c_{1,1}(\delta) -
\lambda_1(\delta) ) c_{3,2}(\delta) - c_{1,2}(\delta)
c_{3,1}(\delta)}{( c_{1,1}(\delta) -
\lambda_1(\delta) ) ( c_{3,3}(\delta) -
\lambda_1(\delta) ) - c_{1,3}(\delta) c_{3,1}(\delta)}
  .
\end{eqnarray*}
Expanding the right-hand sides now furnishes
\[
w_1(\delta)   =   \frac{1}{\sqrt{6}}
\pmatrix{  1 \cr -2 \cr 1
} +
\frac{\sqrt{6}}{24} \cdot\frac{\mathcal{R}_{3,2}^m(x)}{\mathcal
{R}_{3,3}^m(x)}
\pmatrix{  1 \cr 0 \cr -1
}
\cdot\delta+
\pmatrix{  w_{1,1,2} \cr 0 \cr w_{1,3,2}
}
\cdot\delta^2 + O( \delta^3 )   ,
\]
with
\[
w_{1,1,2} + w_{1,3,2}   =
\frac{1}{4\sqrt{6}} \cdot
\frac{\mathcal{R}_{3,1}^m(x)}{\mathcal{R}_{3,3}^m(x)}   .
\]
This finally implies
\[
\tau(\delta) \cdot w_1(\delta)   =
\frac{\mathcal{R}_{3,1}^m(x) \mu(x) - \mathcal{R}_{3,2}^m(x) \mu
^\prime(x) +
\mathcal{R}_{3,3}^m(x)\mu^{\prime\prime}(x)}{4\sqrt{6}
\mathcal{R}_{3,3}^m(x)}
\cdot\delta^2 +
O( \delta^3 )
  ,
\]
and together with
\[
\tau(\delta) \cdot v_1(\delta) =
\frac{-\sqrt{6}   v_{1,2}(\delta)}{2} \cdot
\tau(\delta) \cdot w_1(\delta)
\quad  \mbox{and}\quad
\lim_{\delta\to0} \frac{-\sqrt{6}   v_{1,2}(\delta)}{2} = 1
\]
this establishes the asymptotic behavior of~$\tau(\delta) \cdot
v_1(\delta)$.

Finally, we turn our attention to the second eigenvector. Following our
above approach, we write $v_2(\delta) = (v_{2,1}(\delta),
v_{2,2}(\delta),
v_{2,3}(\delta))^t$, and consider the functions
\[
w_{2,k}(\delta) = \frac{1}{\sqrt{2} \cdot v_{2,1}(\delta)} \cdot
v_{2,k}(\delta)
\qquad  \mbox{for }
k = 1,2,3   .
\]
Then the vector $w_2(\delta) = (w_{2,1}(\delta), w_{2,2}(\delta),
w_{2,3}(\delta))^t$ is defined for sufficiently small $\delta> 0$,
and for these~$\delta$ we have
\[
w_{2,1}(\delta) = \tfrac{1}{\sqrt{2}}
\quad  \mbox{as well as}\quad
\bigl( C(\delta) - \lambda_2(\delta) I \bigr) w_2(\delta) = 0   .
\]
Using again the abbreviation $C(\delta) = (c_{i,j}(\delta))_{i,j=1,2,3}$,
the latter system is equivalent to
\begin{eqnarray*}
\bigl( c_{2,2}(\delta) - \lambda_2(\delta) \bigr) w_{2,2}(\delta
) +
c_{2,3}(\delta) w_{2,3}(\delta) & = & -\tfrac{1}{\sqrt{2}}
\cdot c_{2,1}(\delta)   , \\
c_{3,2}(\delta) w_{2,2}(\delta) + \bigl( c_{3,3}(\delta) -
\lambda_2(\delta) \bigr) w_{2,3}(\delta) & = & -\tfrac{1}{\sqrt{2}}
\cdot c_{3,1}(\delta)   ,
\end{eqnarray*}
which immediately implies
\begin{eqnarray*}
w_{2,2}(\delta) & = &
-\frac{1}{\sqrt{2}} \cdot\frac{( c_{3,3}(\delta) -
\lambda_2(\delta) ) c_{2,1}(\delta) - c_{2,3}(\delta)
c_{3,1}(\delta)}{( c_{2,2}(\delta) -
\lambda_2(\delta) ) ( c_{3,3}(\delta) -
\lambda_2(\delta) ) - c_{2,3}(\delta) c_{3,2}(\delta)}
  , \\
w_{2,3}(\delta) & = &
-\frac{1}{\sqrt{2}} \cdot\frac{( c_{2,2}(\delta) -
\lambda_2(\delta) ) c_{3,1}(\delta) - c_{3,2}(\delta)
c_{2,1}(\delta)}{( c_{2,2}(\delta) -
\lambda_2(\delta) ) ( c_{3,3}(\delta) -
\lambda_2(\delta) ) - c_{2,3}(\delta) c_{3,2}(\delta)}
  .
\end{eqnarray*}
Expanding the right-hand sides now furnishes
\[
w_2(\delta)   =   \frac{1}{\sqrt{2}}
\pmatrix{  1 \cr 0 \cr -1}+
\pmatrix{ 0 \cr\vspace*{2pt} w_{2,2,1} \vspace*{2pt}\cr w_{2,3,1}
}
\cdot\delta+ O( \delta^2 )
\]
with
\[
w_{2,2,1} + w_{2,3,1}   =
\frac{1}{\sqrt{2}} \cdot
\frac{R_{1,0}(x)}{R_{0,0}(x)}   .
\]
This finally implies
\[
\tau(\delta) \cdot w_2(\delta)   =
\frac{R_{1,0}(x) \mu(x) - R_{0,0}(x) \mu^\prime(x)}{\sqrt{2}
R_{0,0}(x)}
\cdot\delta+
O( \delta^2 )
\]
and together with
\[
\tau(\delta) \cdot v_2(\delta) =
\sqrt{2}   v_{2,1}(\delta) \cdot\tau(\delta) \cdot w_2(\delta)
 \quad \mbox{and}\quad
\lim_{\delta\to0} \sqrt{2}   v_{2,1}(\delta) = 1
\]
this establishes the asymptotic behavior of~$\tau(\delta) \cdot
v_2(\delta)$.
\end{pf}

After these preparations, we are finally in a position to prove our
second main result. As mentioned in Section~\ref{secintro}, this result
provides a general means for determining the location of sampling points
of random fields in such a way that the topology of the underlying nodal
sets is correctly recognized with the largest probability. In addition,
the sampling density can readily be determined from derivatives of the
spatial correlation function of the random field.

\begin{pf*}{Proof of Theorem~\ref{thmCorrSamp}}
Due to our assumptions, the random variable $u(x,\cdot) \dvtx  \Omega
\to\mathbb{R}$ is normally distributed with mean~$0$ and its
variance~$R_{0,0}(x)$
is positive for each $x \in[a,b]$ due to (G2).
This immediately implies~(A1). Furthermore, (A2) follows
readily from~\cite{adler81a}, Theorem~3.2.1. Thus, in order to apply
Theorem~\ref{thmTopSamp} we only have to verify~(A3).

For this, we apply Corollary~\ref{corPToolA} with $n = 3$
and sign vector $(s_1,s_2,s_3) = (1,-1,1)$. Fix $x \in[a,b)$ and
consider the $\delta$-dependent three-dimensional random
vector~$T(\delta)$ defined in~(\ref{lemCorrSamp1}).
Then according to Lemma~\ref{lemCorrSamp}, this random vector
satisfies all of the assumptions of Proposition~\ref{propPToolA}
and Corollary~\ref{corPToolA} with
\[
\det C(\delta) = \frac{1}{64} \cdot\det\mathcal{R}(x) \cdot\delta
^6 +
O( \delta^7 )
 \quad \mbox{and}\quad
\lambda_1(\delta) = \frac{\det\mathcal{R}(x)}{96   \mathcal
{R}_{3,3}^m(x)}
\cdot\delta^4 + O( \delta^5 )
\]
as well as
\begin{eqnarray*}
\alpha_1 & = & \frac{\mathcal{R}_{3,1}^m(x) \mu(x) - \mathcal
{R}_{3,2}^m(x) \mu
^\prime(x)
+ \mathcal{R}_{3,3}^m(x) \mu^{\prime\prime}(x)}{\mathcal
{R}_{3,3}^m(x)^{1/2}
\det\mathcal{R}(x)^{1/2}}   , \\
\alpha_2 & = & \frac{R_{1,0}(x) \mu(x) - R_{0,0}(x) \mu^\prime(x)}
{R_{0,0}(x)^{1/2}   \mathcal{R}_{3,3}^m(x)^{1/2}}
  , \\
\alpha_3 & = & \frac{\mu(x)}{R_{0,0}(x)^{1/2}}   .
\end{eqnarray*}
Applying Corollary~\ref{corPToolA}, we then obtain
\[
\lim_{\delta\to0} \bigl( p_{+1}(x,\delta) + p_{-1}(x,\delta)
\bigr)
\cdot\sqrt{\frac{\det C(\delta)}{\lambda_1(\delta)^3}}   =
\frac{3 \sqrt{6}}{4 \pi} \cdot( 1 + \alpha_1^2 )
\cdot
e^{-( \alpha_2^2 + \alpha_3^2 ) / 2}   ,
\]
where we used the formula for~$S_\alpha^{\pm}$ given
in~(\ref{corPToolA2}). In combination with the above expansions
for~$\det C(\delta)$ and~$\lambda_1(\delta)$, this limit furnishes
\[
p_{+1}(x,\delta) + p_{-1}(x,\delta) =
\frac{( 1 + \alpha_1^2 ) \cdot
e^{-( \alpha_2^2 + \alpha_3^2 ) / 2}}{64 \pi} \cdot
\frac{\det\mathcal{R}(x)}{\mathcal{R}_{3,3}^m(x)^{3/2}} \cdot
\delta^3
+ O( \delta^4 )   .
\]
Thus, assumption~(A3) is satisfied with~$\mathcal{C}_0(x) = 3 \mathcal
{C}(x) / 4$, and
Theorem~\ref{thmCorrSamp} follows now immediately from Theorem~\ref
{thmTopSamp}.
\end{pf*}

\section{Concluding remarks}
At first glance, the title of this paper may appear somewhat misleading or
more ambitious than the results delivered. After all, the techniques of proof
are based on classical probabilistic arguments. However, the results
are new
and the examples of Section~\ref{secappl} demonstrate that they have
interesting
nonintuitive implications.

A reasonable question is why were these results not discovered
sooner. We believe that the answer comes from the fact that we are
approaching the
problem of optimal sampling from the point of view of trying to obtain
topological
information. This point of view had been taken previously in the work of
Adler and Taylor~\cite{adler81a,adlertaylor07a}. Their main focus, however,
was the estimation of excursion probabilities, that is, the likelihood that
a given random function exceeds a certain threshold. In~\cite{adler81a,adlertaylor07a}, it is shown that such excursion probabilities can be
well-approximated by studying the geometry of random sub- or
super-level sets
of random fields. More precisely, it is shown that the expected value
of the
Euler characteristic of super-level sets approximates excursion probabilities
for large values of the threshold, and that it is possible to derive explicit
formulas for the expected values of the Euler characteristic and other intrinsic
volumes of nodal domains of random fields.

All of the above results concern the intrinsic volumes of the nodal
domains---which are additive set functionals, and therefore computable
via local considerations alone~\cite{klainrota97a,santalo04a}. In contrast,
in previous work~\cite{gameiroetal05a} we have demonstrated that the
homological
analysis of patterns of nodal sets can uncover phenomena that cannot be captured
using for example only the Euler characteristic. The~more detailed information
on the geometry of patterns encoded in homology is an inherently global quantity
and cannot be computed through local considerations alone. On the other hand,
recent computational advances allow for the fast computation of homological
information based on discretized nodal domains. For this reason, we
focus on
the interface between the discretization and the underlying nodal
domain, rather
than the homology of the nodal domain directly, and then quantify the likelihood
of error in the probabilistic setting. In this sense, our approach complements
the above-mentioned results on the geometry of random fields by Adler and
Taylor~\cite{adler81a,adlertaylor07a}.

Given the current activity surrounding the ideas of using topological methods
for data analysis and remote sensing~\cite{desilvaghrist07a,desilvaghrist07b,ghrist08a}, we believe the importance of this perspective will grow.
Thus, the
title of our paper is chosen in part to encourage the interested reader
to consider the natural generalizations of this work to higher-dimensional
domains where the question becomes one of optimizing the homology of
the generalized nodal sets in terms of homology computed using a complex
derived from a nonuniform sampling of space.

%

\printaddresses


\begin{thebibliography}{30}

\bibitem{adler81a}
%
\begin{bbook}[mr]
\bauthor{\bsnm{Adler},~\bfnm{Robert~J.}\binits{R.~J.}}
(\byear{1981}).
\btitle{The~Geometry of Random Fields}.
\bpublisher{Wiley}, \baddress{Chichester}.
\bid{mr={611857}}
\end{bbook}
%
\endbibitem

\bibitem{adlertaylor07a}
%
\begin{bbook}[mr]
\bauthor{\bsnm{Adler},~\bfnm{Robert~J.}\binits{R.~J.}} \AND
\bauthor{\bsnm{Taylor},~\bfnm{Jonathan~E.}\binits{J.~E.}}
(\byear{2007}).
\btitle{Random Fields and Geometry}.
\bpublisher{Springer}, \baddress{New York}.
\bid{mr={2319516}}
\end{bbook}
%
\endbibitem

\bibitem{bauer96a}
%
\begin{bbook}[mr]
\bauthor{\bsnm{Bauer},~\bfnm{Heinz}\binits{H.}}
(\byear{1996}).
\btitle{Probability Theory}.
\bseries{de Gruyter Studies in Mathematics}
\bvolume{23}.
\bpublisher{de Gruyter}, \baddress{Berlin}.
\bid{mr={1385460}}
\end{bbook}
%
\endbibitem

\bibitem{bharuchareids86a}
%
\begin{bbook}[mr]
\bauthor{\bsnm{Bharucha-Reid},~\bfnm{A.~T.}\binits{A.~T.}} \AND
\bauthor{\bsnm{Sambandham},~\bfnm{M.}\binits{M.}}
(\byear{1986}).
\btitle{Random Polynomials}.
\bpublisher{Academic Press}, \baddress{Orlando, FL}.
\bid{mr={856019}}
\end{bbook}
%
\endbibitem

\bibitem{couranthilbert53a}
%
\begin{bbook}[mr]
\bauthor{\bsnm{Courant},~\bfnm{R.}\binits{R.}} \AND
\bauthor{\bsnm{Hilbert},~\bfnm{D.}\binits{D.}}
(\byear{1953}).
\btitle{Methods of Mathematical Physics. {V}ol. {I}}.
\bpublisher{Interscience Publishers}, \baddress{New York}.
\bid{mr={0065391}}
\end{bbook}
%
\endbibitem

\bibitem{cramerleadbetter04a}
%
\begin{bbook}[mr]
\bauthor{\bsnm{Cram{\'e}r},~\bfnm{Harald}\binits{H.}} \AND
\bauthor{\bsnm{Leadbetter},~\bfnm{M.~R.}\binits{M.~R.}}
(\byear{2004}).
\btitle{Stationary and Related Stochastic Processes}.
\bpublisher{Dover}, \baddress{Mineola, NY}.
\bnote{Reprint of the 1967 original}.
\bid{mr={2108670}}
\end{bbook}
%
\endbibitem

\bibitem{dayetal07a}
%
\begin{barticle}[mr]
\bauthor{\bsnm{Day},~\bfnm{Sarah}\binits{S.}},
\bauthor{\bsnm{Kalies},~\bfnm{William~D.}\binits{W.~D.}},
\bauthor{\bsnm{Mischaikow},~\bfnm{Konstantin}\binits{K.}} \AND
\bauthor{\bsnm{Wanner},~\bfnm{Thomas}\binits{T.}}
(\byear{2007}).
\btitle{Probabilistic and numerical validation of homology
computations for
nodal domains}.
\bjournal{Electron. Res. Announc. Amer. Math. Soc.}
\bvolume{13}
\bpages{60--73 (electronic)}.
\bid{doi={10.1090/S1079-6762-07-00175-8}, mr={2320683}}
\end{barticle}
%
\endbibitem

\bibitem{desilvaghrist07a}
%
\begin{barticle}[mr]
\bauthor{\bparticle{de~}\bsnm{Silva},~\bfnm{Vin}\binits{V.}} \AND
\bauthor{\bsnm{Ghrist},~\bfnm{Robert}\binits{R.}}
(\byear{2007}).
\btitle{Coverage in sensor networks via persistent homology}.
\bjournal{Algebr. Geom. Topol.}
\bvolume{7}
\bpages{339--358}.
\bid{doi={10.2140/agt.2007.7.339}, mr={2308949}}
\end{barticle}
%
\endbibitem

\bibitem{desilvaghrist07b}
%
\begin{barticle}[mr]
\bauthor{\bparticle{de~}\bsnm{Silva},~\bfnm{Vin}\binits{V.}} \AND
\bauthor{\bsnm{Ghrist},~\bfnm{Robert}\binits{R.}}
(\byear{2007}).
\btitle{Homological sensor networks}.
\bjournal{Notices Amer. Math. Soc.}
\bvolume{54}
\bpages{10--17}.
\bid{mr={2275921}}
\end{barticle}
%
\endbibitem

\bibitem{dunnage66a}
%
\begin{barticle}[mr]
\bauthor{\bsnm{Dunnage},~\bfnm{J.~E.~A.}\binits{J.~E.~A.}}
(\byear{1966}).
\btitle{The~number of real zeros of a random trigonometric polynomial}.
\bjournal{Proc. London Math. Soc. (3)}
\bvolume{16}
\bpages{53--84}.
\bid{mr={0192532}}
\end{barticle}
%
\endbibitem

\bibitem{edelmankostlan95a}
%
\begin{barticle}[mr]
\bauthor{\bsnm{Edelman},~\bfnm{Alan}\binits{A.}} \AND
\bauthor{\bsnm{Kostlan},~\bfnm{Eric}\binits{E.}}
(\byear{1995}).
\btitle{How many zeros of a random polynomial are real?}
\bjournal{Bull. Amer. Math. Soc. (N.S.)}
\bvolume{32}
\bpages{1--37}.
\bid{doi={10.1090/S0273-0979-1995-00571-9}, mr={1290398}}
\end{barticle}
%
\endbibitem

\bibitem{farahmand98a}
%
\begin{bbook}[mr]
\bauthor{\bsnm{Farahmand},~\bfnm{Kambiz}\binits{K.}}
(\byear{1998}).
\btitle{Topics in Random Polynomials}.
\bseries{Pitman Research Notes in Mathematics Series}
\bvolume{393}.
\bpublisher{Longman}, \baddress{Harlow}.
\bid{mr={1679392}}
\end{bbook}
%
\endbibitem

\bibitem{gameiroetal04a}
%
\begin{barticle}[mr]
\bauthor{\bsnm{Gameiro},~\bfnm{Marcio}\binits{M.}},
\bauthor{\bsnm{Mischaikow},~\bfnm{Konstantin}\binits{K.}} \AND
\bauthor{\bsnm{Kalies},~\bfnm{William}\binits{W.}}
(\byear{2004}).
\btitle{Topological characterization of spatial-temporal chaos}.
\bjournal{Phys. Rev. E (3)}
\bvolume{70}
\bpages{035203, 4}.
\bid{doi={10.1103/PhysRevE.70.035203}, mr={2129999}}
\end{barticle}
%
\endbibitem

\bibitem{gameiroetal05a}
%
\begin{barticle}[auto:SpringerTagBib|2009-01-14|16:51:27]
\bauthor{\bsnm{Gameiro},~\bfnm{M.}\binits{M.}},
\bauthor{\bsnm{Mischaikow},~\bfnm{K.}\binits{K.}} \AND
\bauthor{\bsnm{Wanner},~\bfnm{T.}\binits{T.}}
(\byear{2005}).
\btitle{Evolution of pattern complexity in the {C}ahn--{H}illiard
theory of
phase separation}.
\bjournal{Acta Materialia}
\bvolume{53}
\bpages{693--704}.
\bid{doi={10.1016/j.actamat.2004.10.022}}
\end{barticle}
%
\endbibitem

\bibitem{ghrist08a}
%
\begin{barticle}[mr]
\bauthor{\bsnm{Ghrist},~\bfnm{Robert}\binits{R.}}
(\byear{2008}).
\btitle{Barcodes: The~persistent topology of data}.
\bjournal{Bull. Amer. Math. Soc. (N.S.)}
\bvolume{45}
\bpages{61--75 (electronic)}.
\bid{doi={10.1090/S0273-0979-07-01191-3}, mr={2358377}}
\end{barticle}
%
\endbibitem

\bibitem{kac43a}
%
\begin{barticle}[mr]
\bauthor{\bsnm{Kac},~\bfnm{M.}\binits{M.}}
(\byear{1943}).
\btitle{On the average number of real roots of a random algebraic equation}.
\bjournal{Bull. Amer. Math. Soc.}
\bvolume{49}
\bpages{314--320}.
\bid{mr={0007812}}
\end{barticle}
%
\endbibitem

\bibitem{kac49a}
%
\begin{barticle}[mr]
\bauthor{\bsnm{Kac},~\bfnm{M.}\binits{M.}}
(\byear{1949}).
\btitle{On the average number of real roots of a random algebraic equation.
{II}}.
\bjournal{Proc. London Math. Soc. (2)}
\bvolume{50}
\bpages{390--408}.
\bid{mr={0030713}}
\end{barticle}
%
\endbibitem

\bibitem{kaczynskietal04a}
%
\begin{bbook}[mr]
\bauthor{\bsnm{Kaczynski},~\bfnm{Tomasz}\binits{T.}},
\bauthor{\bsnm{Mischaikow},~\bfnm{Konstantin}\binits{K.}} \AND
\bauthor{\bsnm{Mrozek},~\bfnm{Marian}\binits{M.}}
(\byear{2004}).
\btitle{Computational Homology}.
\bseries{Applied Mathematical Sciences}
\bvolume{157}.
\bpublisher{Springer}, \baddress{New York}.
\bid{mr={2028588}}
\end{bbook}
%
\endbibitem

\bibitem{kahane85a}
%
\begin{bbook}[mr]
\bauthor{\bsnm{Kahane},~\bfnm{Jean-Pierre}\binits{J.-P.}}
(\byear{1985}).
\btitle{Some Random Series of Functions},
\bedition{2nd} ed.
\bseries{Cambridge Studies in Advanced Mathematics}
\bvolume{5}.
\bpublisher{Cambridge Univ. Press}, \baddress{Cambridge}.
\bid{mr={833073}}
\end{bbook}
%
\endbibitem

\bibitem{klainrota97a}
%
\begin{bbook}[mr]
\bauthor{\bsnm{Klain},~\bfnm{Daniel~A.}\binits{D.~A.}} \AND
\bauthor{\bsnm{Rota},~\bfnm{Gian-Carlo}\binits{G.-C.}}
(\byear{1997}).
\btitle{Introduction to Geometric Probability}.
\bpublisher{Cambridge Univ. Press}, \baddress{Cambridge}.
\bid{mr={1608265}}
\end{bbook}
%
\endbibitem

\bibitem{krishanetal07a}
%
\begin{barticle}[auto:SpringerTagBib|2009-01-14|16:51:27]
\bauthor{\bsnm{Krishan},~\bfnm{K.}\binits{K.}},
\bauthor{\bsnm{Gameiro},~\bfnm{M.}\binits{M.}},
\bauthor{\bsnm{Mischaikow},~\bfnm{K.}\binits{K.}},
\bauthor{\bsnm{Schatz},~\bfnm{M.}\binits{M.}},
\bauthor{\bsnm{Kurtuldu},~\bfnm{H.}\binits{H.}} \AND
\bauthor{\bsnm{Madruga},~\bfnm{S.}\binits{S.}}
(\byear{2007}).
\btitle{Homology and symmetry breaking in {R}ayleigh--{B}enard convection:
{E}xperiments and simulations}.
\bjournal{Phys. Fluids}
\bvolume{19}
\bpages{117105}.
\end{barticle}
%
\endbibitem

\bibitem{marcuspisier81a}
%
\begin{bbook}[mr]
\bauthor{\bsnm{Marcus},~\bfnm{Michael~B.}\binits{M.~B.}} \AND
\bauthor{\bsnm{Pisier},~\bfnm{Gilles}\binits{G.}}
(\byear{1981}).
\btitle{Random {F}ourier Series with Applications to Harmonic Analysis}.
\bseries{Annals of Mathematics Studies}
\bvolume{101}.
\bpublisher{Princeton Univ. Press}, \baddress{Princeton, NJ}.
\bid{mr={630532}}
\end{bbook}
%
\endbibitem

\bibitem{mischaikowwanner07a}
%
\begin{barticle}[mr]
\bauthor{\bsnm{Mischaikow},~\bfnm{Konstantin}\binits{K.}} \AND
\bauthor{\bsnm{Wanner},~\bfnm{Thomas}\binits{T.}}
(\byear{2007}).
\btitle{Probabilistic validation of homology computations for nodal domains}.
\bjournal{Ann. Appl. Probab.}
\bvolume{17}
\bpages{980--1018}.
\bid{doi={10.1214/105051607000000050}, mr={2326238}}
\end{barticle}
%
\endbibitem

\bibitem{rozanov98a}
%
\begin{bbook}[mr]
\bauthor{\bsnm{Rozanov},~\bfnm{Yu.~A.}\binits{Y.~A.}}
(\byear{1998}).
\btitle{Random Fields and Stochastic Partial Differential Equations}.
\bseries{Mathematics and Its Applications}
\bvolume{438}.
\bpublisher{Kluwer Academic}, \baddress{Dordrecht}.
\bid{mr={1629699}}
\end{bbook}
%
\endbibitem

\bibitem{santalo04a}
%
\begin{bbook}[mr]
\bauthor{\bsnm{Santal{\'o}},~\bfnm{Luis~A.}\binits{L.~A.}}
(\byear{2004}).
\btitle{Integral Geometry and Geometric Probability},
\bedition{2nd} ed.
\bpublisher{Cambridge Univ. Press}, \baddress{Cambridge}.
\bid{mr={2162874}}
\end{bbook}
%
\endbibitem

\bibitem{sidorovetal02a}
%
\begin{bbook}[mr]
\bauthor{\bsnm{Sidorov},~\bfnm{Nikolay}\binits{N.}},
\bauthor{\bsnm{Loginov},~\bfnm{Boris}\binits{B.}},
\bauthor{\bsnm{Sinitsyn},~\bfnm{Aleksandr}\binits{A.}} \AND
\bauthor{\bsnm{Falaleev},~\bfnm{Michail}\binits{M.}}
(\byear{2002}).
\btitle{Lyapunov--{S}chmidt Methods in Nonlinear Analysis and Applications}.
\bseries{Mathematics and Its Applications}
\bvolume{550}.
\bpublisher{Kluwer Academic}, \baddress{Dordrecht}.
\bid{mr={1959647}}
\end{bbook}
%
\endbibitem

\bibitem{torquato02b}
%
\begin{bbook}[mr]
\bauthor{\bsnm{Torquato},~\bfnm{Salvatore}\binits{S.}}
(\byear{2002}).
\btitle{Random Heterogeneous Materials: Microstructure and Macroscopic Properties}.
\bseries{Interdisciplinary Applied Mathematics}
\bvolume{16}.
\bpublisher{Springer}, \baddress{New York}.
\bid{mr={1862782}}
\end{bbook}
%
\endbibitem

\bibitem{vainbergtrenogin74a}
%
\begin{bbook}[mr]
\bauthor{\bsnm{Va{\u\i}nberg},~\bfnm{M.~M.}\binits{M.~M.}} \AND
\bauthor{\bsnm{Trenogin},~\bfnm{V.~A.}\binits{V.~A.}}
(\byear{1974}).
\btitle{Theory of Branching of Solutions of Non-linear Equations}.
\bpublisher{Noordhoff International Publishing}, \baddress{Leyden}.
\bid{mr={0344960}}
\end{bbook}
%
\endbibitem

\bibitem{vanmarcke83a}
%
\begin{bbook}[mr]
\bauthor{\bsnm{Vanmarcke},~\bfnm{Erik}\binits{E.}}
(\byear{1983}).
\btitle{Random Fields}.
\bpublisher{MIT Press}, \baddress{Cambridge, MA}.
\bid{mr={761904}}
\end{bbook}
%
\endbibitem

\bibitem{wilkinson88a}
%
\begin{bbook}[mr]
\bauthor{\bsnm{Wilkinson},~\bfnm{J.~H.}\binits{J.~H.}}
(\byear{1988}).
\btitle{The~Algebraic Eigenvalue Problem}.
\bpublisher{Oxford Univ. Press}, \baddress{New York}.
\bid{mr={950175}}
\end{bbook}
%
\endbibitem

\end{thebibliography}
\end{document}